\documentclass[amsthm]{elsart}

\usepackage{yjsco}
\usepackage{natbib}

\usepackage{multirow} % for table multirow
\usepackage{enumitem} % for parsep, itemsep...
\usepackage{xspace} % for \xspace

\usepackage{amsmath,amssymb}
\usepackage{amsthm}
\usepackage{bm}
\usepackage{pdfsync}
\usepackage{subfigure}
\usepackage[dvipsnames]{xcolor}
\usepackage[english]{babel}

\usepackage[T1]{fontenc}
\usepackage{yfonts}

\theoremstyle{plain}
\newtheorem{theorem}{Theorem}[section]
\newtheorem{proposition}[theorem]{Proposition}
\newtheorem{lemma}[theorem]{Lemma}
\newtheorem{corollary}[theorem]{Corollary}

\theoremstyle{definition}
\newtheorem{example}[theorem]{Example}
\newtheorem{definition}[theorem]{Definition}
\newtheorem{remark}[theorem]{Remark}
\newtheorem{algo}[theorem]{Algorithm}

\newenvironment{ABPRAlgorithm}
  {\goodbreak\bigskip\hrule\begin{algo}}
  {\end{algo}\hrule\medskip}

\newcommand{\minpoly}[1]{{\mu_{#1}}}
\newcommand{\charpoly}[1]{{\chi_{#1}}}
\let\phi=\varphi
\let\rho=\varrho
\let\theta=\vartheta
\let\epsilon=\varepsilon

\def\LT{\mathop{\rm LT}\nolimits}

\def\NF{\mathop{\rm NF}\nolimits}

\def\Supp{\mathop{\rm Supp}\nolimits}

\def\den{\mathop{\rm den}\nolimits}
\def\mod{\mathop{\rm mod}\,}
\def\lcm{\mathop{\rm lcm}\nolimits}
\def\gcd{\mathop{\rm gcd}\nolimits}
\def\notdiv{{\not|}\,}
\newcommand \ie {\textit{i.e.}}
\newcommand \eg {\textit{e.g.}}
\newcommand \ideal[1] {\langle #1 \rangle}

\newcommand \TotalSplit{\texttt{TotalSplit}}
\newcommand \PartialSplit{\texttt{PartialSplit}}
\def\sqfree{\mathop{\rm sqfree}\nolimits}
\def\Rad{\mathop{\rm Rad}\nolimits}

\def\Frob{\mathop{\rm Frob}\nolimits}

\def\phi{{\varphi}}

\newcommand \FF {{\mathbb F}}
\newcommand \NN {{\mathbb N}}
\newcommand \QQ {{\mathbb Q}}

\newcommand \TT {{\mathbb T}}
\newcommand \ZZ {{\mathbb Z}}

\def\To{\longrightarrow}
\def\TTo#1{\mathop{\longrightarrow}\limits^{#1}}

\let\To=\longrightarrow
\def\TTo#1{\mathop{\longrightarrow}\limits ^{#1}}

\def\calc{_{\rm calc}}
\def\crt{_{\rm crt}}

\def\longiso{\,\smash{\TTo{\lower 7pt\hbox{$\scriptstyle\sim$}}}\,}

\def\tfrac #1#2{{\textstyle\frac{#1}{#2}}}

\def\cocoa{\mbox{\rm
C\kern-.13em o\kern-.07 em C\kern-.13em o\kern-.15em A}\xspace}
\def\apcocoa{\mbox{\rm
A\kern-0.13em p\kern -0.07em C\kern-.13em o\kern-.07 em C\kern-.13em
o\kern-.15em A}}

\begin{document}

\begin{frontmatter}

\title{Computing and Using Minimal Polynomials}

\author{John Abbott}
\address{\scriptsize Institut f\"ur Mathematik, \ Universit\"at
  Kassel, Germany}
\ead{abbott@dima.unige.it}

\author{Anna Maria Bigatti}
\address{\scriptsize Dip. di Matematica,
\ Universit\`a degli Studi di Genova, \ Via
Dodecaneso 35,\
I-16146\ Genova, Italy}
\ead{bigatti@dima.unige.it}

\author{Elisa Palezzato}
\address{\scriptsize Dip. di Matematica,
\ Universit\`a degli Studi di Genova, \ Via
Dodecaneso 35,\
I-16146\ Genova, Italy}
\ead{bigatti@dima.unige.it}

\author{Lorenzo Robbiano}
\address{\scriptsize Dip. di Matematica,
\ Universit\`a degli Studi di Genova, \ Via
Dodecaneso 35,\
I-16146\ Genova, Italy}
\ead{robbiano@dima.unige.it}

\thanks{This research was partly supported by
the project H2020-FETOPN-2015-CSA\_712689 of the European Union.\quad
J.~Abbott is an INdAM-COFUND Marie Curie Fellow.\\
Published on JSC \texttt{https://doi.org/10.1016/j.jsc.2019.07.022}
\\
\copyright 2019. This manuscript version is made available under the CC-BY-NC-ND 4.0 license \texttt{http://creativecommons.org/licenses/by-nc-nd/4.0/}
}
%the ``National Group for Algebraic and Geometric Structures, and
%their Applications''(GNSAGA -- INDAM)

\bigskip
\begin{abstract}
Given a zero-dimensional ideal $I$ in a polynomial ring,
many computations start by finding univariate polynomials in $I$.
Searching for a univariate polynomial in $I$ is a particular case of
considering the minimal polynomial of an element in $P/I$.
It is well known that minimal polynomials  may be computed via elimination,
therefore this is considered to be a ``resolved problem''.
But being the key of so many computations, it is worth investigating
its meaning, its optimization, its applications
(\eg~testing if a zero-dimensional ideal is radical, primary or maximal).
We present efficient algorithms for computing the minimal
  polynomial of an element of $P/I$.
For the specific case
  where the coefficients are in $\QQ$, we show how to use modular
  methods to obtain a guaranteed result.  
We also present some
  applications of minimal polynomials, namely algorithms for computing
  radicals and primary decompositions of zero-dimensional ideals, and
  also for testing radicality and maximality.
\end{abstract}

\date{\today}


\end{frontmatter}

%\tableofcontents

\goodbreak
%%%%%%%%%%%%%%%%%%%%%%%%%%%%%%%%
\section{Introduction}
\label{MinPoly:sec:intro}
%%%%%%%%%%%%%%%%%%%%%%%%%%%%%%%%

This paper describes both theoretical and practical aspects of computing minimal polyomials, with particular emphasis on our practical
implementation in CoCoALib \citep{cocoalib} and \cocoa \citep{cocoa5}.
It is organized into three parts:
computing minimal polynomials over~$\FF_p$,
computing minimal polynomials over~$\QQ$, and
using minimal polynomials.

\smallskip

In linear algebra it is frequently necessary to use non-linear objects such as
minimal and characteristic polynomials since they capture fundamental information
about endomorphisms of finite-dimensional vector spaces.
It is well-known that if $K$ is a field and $R$ is a zero-dimensional affine $K$-algebra
(\ie~a zero-dimensional algebra of the form $R=K[x_1, \dots, x_n]/I$)
then $R$ is a finite-dimensional $K$-vector space~---~\eg~see Proposition~3.7.1 of~\citep{KR1}.
Consequently, it is not surprising that minimal and characteristic polynomials  can be
successfully used to detect properties of~$R$.

This point of view was taken systematically
in the book by~\citet{KR3}, where the particular importance of minimal polynomials
(rather greater than that of characteristic polynomials) emerged quite clearly.
The book clarified one main advantage of minimal polynomials over characteristic
polynomials, namely that minimal polynomials generalise to families of pairwise
commuting endomorphisms, while characteristic polynomials do not.
The book also described several algorithms which use minimal
polynomials as a crucial tool.  The approach taken there
was a good source of inspiration for our research, so
we decided to delve into the theory of minimal polynomials,
their uses, and their applications.

% \cancel{This was the first step, and immediately after we became aware of the problem of
% getting more source of inspiration from the existing literature.}
Being such fundamental tools in linear algebra,
minimal polynomials have played a prominent role
in several branches of research over the last two centuries.
As a consequence, it is impractical
to track down all contributions to this theory.
Two further sources of particular inspiration were the paper by~\citet{La}, where minimal polynomials
were systematically used in the process of solving zero-dimensional systems of polynomial equations,
and the paper by~\citet{BSS}, which contains a fine analysis of the complexity of computing minimal polynomials.

\subsubsection*{Relevance to the {\sf SC$^2$} community}

It is not immediately apparent how minimal polynomials could be relevant to
  the \textsf{SC$^2$} community.
However, our methods for computing them efficiently are very helpful for
solving polynomial systems, 
%and this is a 
and are basic building blocks also for a number of other
applications in the \textsf{SC$^2$} community,
as introduced in~\citep{AbbBig2017}.

%Another
For instance, one
 fundamental tool in the algebraic approach to tackling \textsf{SC$^2$} problems is
 \textit{Cylindrical Algebraic Decomposition}, often abbreviated to \textit{CAD} (first introduced in~\citep{collins-cad}).
The software CArL/SMT-RAT by~\citet{KA2018} employs Lazard's variant of CAD
\citep{lazard-CAD}
which requires polynomial factorization over algebraic field extensions.
We use our efficient algorithms for minimal polynomials to compute quickly
the primary decomposition of zero-dimensional ideals (see
Section~\ref{MinPoly:sec:PrimDec} and
Table~\ref{MinPoly:tab:SingularCoCoA}), which in turn give a good method
for factorizing polynomials over algebraic field extensions.
We presented the relevant theory and our CoCoALib implementation
in~\citep{AbbBigPal2018}; further background material is in Section~5.3.B of~\citep{KR3},
  Indeed, our CoCoALib implementation is used by CArL/SMT-RAT.

% \renzo{The theoretical background of this technique, 
% together with suitable examples, is
% explained in~\citet{KR3}, Section 5.3.B. As said, the main ingredient is the 
% computation of the primary decomposition of a suitable zero-dimensional ideal
% which we do according to the
% detailed explanation in Subsection~\ref{MinPoly:sec:PrimDec}, and whose excellent
% performance is confirmed by Table~\ref{MinPoly:tab:SingularCoCoA}).
% }

\subsubsection*{Fundamentals of our approach}

The first step of our approach was to devise and implement good algorithms
for computing the minimal polynomial of an element of $R$
and that of a $K$-endomorphism of~$R$ (see Algorithms~\ref{MinPoly:alg:MinPolyQuotMat}, and~\ref{MinPoly:alg:MinPolyQuotDef}).
They are described in Section~\ref{MinPoly:sec:algms}, and refine similar algorithms examined in the book by~\citet{KR3}.
They have been implemented in CoCoALib~\citep{cocoalib}, and are accessible from \cocoa~\citep{cocoa5}, as
indeed are all other algorithms described in this paper.
Although the theoretical content of Section~\ref{MinPoly:sec:algms} is
essentially elementary,
many remarks on implementation details are given to get
the good performance shown in Table~\ref{MinPoly:tab:PrimeFields}.
%(see for instance Remark~\ref{rem:important}).

\subsubsection*{Efficient computation by modular techniques}

%\smallskip
Section~\ref{MinPoly:sec:modularapproach}
constitutes a contribution of great practical significance:
it addresses the problem of computing minimal polynomials of elements of
an affine $\QQ$-algebra using a modular approach.
The technique of using modular reduction has a long tradition,
see for instance these articles by~\citet{Pa}, \citet{Wi}, \citet{Gr}, \citet{NY[17]},
\citet{NY[18]}, \citet{NY2018}, \citet{AoyNor}, \citet{Nor2002}, \citet{IPS} and \citet{M04}.
% \john{Preferisco qui non dire niente al riguardo di ``probably correct''
%   perche' si presuppone che il risultato sia sempre giusto a meno che
%   non venga esplicitamente detto il contrario.  Forse si puo` dire
%   qualcosa tipo: We refer to Remark~\ref{MinPoly:rem:Certified answer} for
%   a brief discussion about a faster variant of our approach which
% produces a ``probably correct'' result.}
%\cancel{
Usually modular methods are associated with results which are correct
\textit{only with high probability}.
In contrast, with the approach we introduce here, the minimal
polynomial is guaranteed to be correct (see Remark~\ref{MinPoly:rem:Certified answer}).
%}

Why might a modular approach be useful for computing minimal polynomials?
Let us have a look at an example.  Let $\mathbb{X} = \{p_1, p_2,p_3,p_4\}$  be
the set of the following four points in $\mathbb{A}^3_\QQ$:
$p_1 = (1,3,0)$, $p_2=(1,0,4)$, $p_3 = (-5,7,1)$, $p_4=(1,0,0)$.
The vanishing ideal of $\mathbb{X}$ in the polynomial ring $P=\QQ[x,y,z]$ is
$I =
\ideal{z^2 - \tfrac{1}{2} x -4 z + \tfrac{1}{2},\allowbreak
  yz + \tfrac{7}{6} x - \tfrac{7}{6},\;
  y^2 + \tfrac{14}{3} x -3 y - \tfrac{14}{3}}$.
We take $f = 2x^2  +3y^4 +5z^6 \in P$.  Then the minimal polynomial of $\bar{f}$ in $P/I$,
\ie~the lowest-degree, monic polynomial $g(z)$ such that $g(f) \in I$,  is
$z^4 -27987\,z^3 +155510626\,z^2 -36732206532\,z +72842594440$.
Even in this ``small'' example, although both the ideal and the element $f$ have simple coefficients,
the minimal polynomial has much larger coefficients.

The modular approach tames the complexity of computing the big coefficients of such
polynomials.  However,
as always happens with a modular approach, various obstacles have to be overcome~---~see for instance the discussion contained in~\citep{ABR}.
In particular, we deal with the notion of reduction of an ideal modulo
$p$, and we do so by introducing the notion of $\sigma$-denominator of
an ideal (see Definition~\ref{MinPoly:def:reductionmodp} and
Theorem~\ref{MinPoly:thm:RGBp}), which enables us to surmount these obstacles.
The reason why we introduce the $\sigma$-denominator is that,
for a given a term ordering,  the
reduced Gr\"obner basis of an ideal is unique, so that the theoretical results
do not depend on the generators of the ideal, but just on the ideal itself.

Hand-in-hand with the modular approach go the notions of \textit{usable, good} and \textit{bad} primes (see Definitions~\ref{MinPoly:def:usableprime}
and \ref{MinPoly:def:goodprime}).  We show that all but finitely many primes are good
(see Theorem~\ref{MinPoly:thm:almostallgood} and Corollary~\ref{MinPoly:cor:badprimes}),
and this paves the way to the construction of the fundamental Algorithm~\ref{MinPoly:alg:modular}.
In the literature, several authors have looked at various notions of bad reduction in similar contexts,
%\cancel{Good and bad primes have several ancestors in the literature},
see for instance the articles by \citet{NY2018},
  \citet{Pa}, \citet{Wi} and~\citet{Ar}.  However, our approach is systematically tied to reduced Gr\"obner
  bases as computationally robust representations of ideals; we study
 more deeply this approach in the preprint~\citep{ABR-Modp}.
The combination of the theoretical results explained in this section with various implementation
details lead to the good practical performance as shown in Table~\ref{MinPoly:tab:rationals}.
For example, the choice of larger ``small primes'' as moduli turned out to be a winning strategy, namely with about $30$ bits (on a 64-bit platform).

\subsubsection*{Uses of minimal polynomials: radical, primary decomposition}

%Section~\ref{MinPoly:sec:timings} presents non-trivial examples of
%minimal polynomials computed with \cocoa, and
Section~\ref{MinPoly:sec:uses} shows how minimal polynomials
can be successfully and efficiently used to compute several important invariants
of zero-dimensional affine $K$-algebras.
More specifically, in Subsection~\ref{MinPoly:subsec:radical}
we describe Algorithms~\ref{MinPoly:alg:IsRadical0Dim}
and \ref{MinPoly:alg:Radical0Dim} which show respectively how to determine whether a zero-dimensional ideal is radical,
and how to compute the radical of a zero-dimensional ideal.
In Subsection~\ref{MinPoly:subsec:IsMaximal} we present several algorithms which determine whether a
zero-dimensional ideal is maximal or primary. The techniques used depend very much
on the field $K$.  The main distinction is between small finite fields and
fields of characteristic zero or big fields of positive characteristic.
In particular, it is noteworthy that in the first case Frobenius spaces
play a fundamental role~---~\eg~see Section~5.2 of the book by~\citet{KR3}.

Finally, in Subsection~\ref{MinPoly:sec:PrimDec} a series of algorithms (see Algorithms~\ref{MinPoly:alg:PDSplitting},
\ref{MinPoly:alg:PDSplittingChar0},  \ref{MinPoly:alg:PrimaryDecompositionCore} and \ref{MinPoly:alg:PrimaryDecomposition0Dim})
describe how to compute the primary decomposition of a zero-dimensional affine $K$-algebra.
They are inspired by the content of Chapter~5 of~\citep{KR3}, but also present many novelties.

% \cancel{ -- moved --
% Our methods for computing efficiently minimal polynomials are very helpful for
% solving polynomial systems, and together with the other
% derived algorithms we believe we can develop tools which will be beneficial to the
% {\sf SC$^2$} community \citet{AbbBig2017}.
% In particular, in the context of Cylindrical Algebraic
% Decomposition~(CAD), in \citet{AbbBigPal2018} we use our primary decomposition
% algorithm for factorizing polynomials over algebraic field extensions.
% }

\subsubsection*{Implementation and timings}

As already mentioned, all the algorithms described in this paper have been implemented in \cocoa. Their
merits are illustrated by the tables of examples contained in Sections~\ref{MinPoly:sec:algms} and~\ref{MinPoly:sec:modularapproach},
and also at the end of Section~\ref{MinPoly:sec:uses}.
The examples were chosen to cover a wide spectrum of
  zero-dimensional affine $K$-algebras;
some are complete intersections, and some are not.
The experiments were performed on a MacBook Pro with Intel Core~i7 processor (clocked at 2.9GHz),
using our implementation in \cocoa~5.

\smallskip
  At the request of the referees we supply comparative timings against two
  (free and open-source) competitors, despite our reservations.
  But is our comparison really fair?  Both Singular (by~\citet{Singular})
  and Macaulay2 (by~\citet{Macaulay2}) are advanced software systems which
  require months to learn thoroughly; for the comparison, we consulted their
  respective manuals, and chose what appear to be the most relevant functions.

We tried to compute our examples with Macaulay2, %(see~\citet{Macaulay2}),
but were unable to compute them in a reasonable time.
In contrast,
we could compute most of the examples with Singular. %(see~\citet{Singular}).
The results are in Table~\ref{MinPoly:tab:SingularCoCoA}, and confirm
the efficiency of our algorithms and implementations (whose outputs, as
already mentioned, are guaranteed to be correct).

%%%%%%%%%%%%%%%%%%%%%%%%%%%%%%%%%%%%%
%\section{Notation, First Definitions, and Examples}
%\label{MinPoly:sec:notation}
\section{Computing Minimal Polynomials}
\label{MinPoly:sec:algms}

Here we introduce the notation and terminology we shall use,
and the definition of minimal polynomial which is the fundamental object
studied in this paper.

Let $K$ be a field, let $P = K[x_1, \dots, x_n]$ be
a polynomial ring in $n$ indeterminates,
and let $\TT^n$ denote the monoid of power-products in $x_1,\ldots,x_n$.
Let $I$ be a  zero-dimensional ideal in~$P$; this implies
that the ring $R=P/I$ is a zero-dimensional affine $K$-algebra,
hence it is a finite dimensional $K$-vector space.
Then, for any $f$ in $P$
there is a linear dependency mod $I$ among the powers of~$f$:
in other words, there is a polynomial $g(z) = \sum_{i=0}^d \lambda_i z^i\in K[z]$
which vanishes modulo $I$ when evaluated at $z=f$.

\begin{example}\label{MinPoly:ex:toy-Min=Char}
Let $P =\FF_{101}[x,y]$,
$I = \ideal{x^2-y,\;  y^2-2x-4}$,
and $f=5x-3y$.
Then %$\minpoly{f,I}(z)$ is
$g(z)= z^4 +18z^2 +48z -23 \in K[z]$
is such that $g(f) = f^4 +18f^2 +48f -23 \in I$,
or equivalently $g(\bar{f}) =0 \in P/I$.
Moreover $g(z)$ is the lowest degree monic polynomial with this property.
\end{example}

\begin{definition}\label{MinPoly:def:minpoly}
  Let $K$ be a field, let $P = K[x_1,\dots, x_n]$, and let $I$ be a
  zero-dimensional ideal.  Given a polynomial $f \in P$, we have a
  $K$-algebra homomorphism $K[z] \to P/I$ given by $z \mapsto f \mod I$.
  The monic generator of the kernel of this homomorphism is called
  the \textbf{minimal polynomial} of $f \mod I$ (or simply ``of $f$'' when the
  ideal $I$ is obvious), and is denoted by
  $\minpoly{f,I}(z)$.
\end{definition}

\begin{remark}\label{rem:minpoly-x}
The particular case of $\minpoly{x_i,I}(x_i)$, where $x_i$ is an
indeterminate, is a very important and popular object when computing:
in fact $\minpoly{x_i,I}(x_i)$ is the lowest degree polynomial in $x_i$
belonging to $I$, that is
$I \cap K[x_i] = \ideal{ \minpoly{x_i,I}(x_i) }$.
% It is well known that this polynomial may be computed via elimination
% of all the other indeterminates $x_j$ (see for example Corollary~3.4.6
% of~\citet{KR1}).  However the algorithm which derives from this
% observation is usually impractically slow.
\end{remark}

\begin{example}\label{MinPoly:ex:toy-Min-neChar}
Let $P =\FF_{101}[x,y]$,
$I = \ideal{y^3 -xy -2y^2 +y,\;  xy^2,\; x^2 -x}$,
and $f=y$.
Then $\minpoly{f,I}(y) = y^4 -2y^3 +y^2 = y^2{\cdot}(y-1)^2$.
% is the lowest degree polynomial in $y$ belonging to $I$.
\end{example}

\begin{example}\label{MinPoly:ex:toy-QBigCoeff}
Let $P =\QQ[x,y]$,
$I = \ideal{x^2-\frac17 y^2-5,\;  y^2+4x-\frac72}$,
%$f_1=y$
and $f = 3x-2y$.
Then
%$\minpoly{f_1,I}(z) = z^4 -\frac{65}7 z^2 -\frac{271}4$ and
$\minpoly{f,I}(z) =
z^4  +\frac{24}{7} z^3  -{\frac{6527}{49}}^{\mathstrut} z^2  +\frac{5868}{7} z
+\frac{10967}{28}$.
We will see in Example~\ref{MinPoly:ex:toy-QBigCoeff-contd}
how we compute it with a modular approach.
% \cancel{,
% showing that rational coefficients can be large, even in
% small examples.}
% \john{In verita` l'esempio usato per giustificare l'uso di un approccio
% modulare e` piu` impressionante.}  VERO - cambiato frase
\end{example}
% \john{Esempio piu` grosso: $\minpoly{f,I}(z) =  z^4 +\tfrac{6912}{4081} z^3 + \tfrac{49024}{4081} z^2 +\tfrac{55296}{4081} z + \tfrac{1458112}{4081}$
%   dove $I =\ideal{x^2 +8 x y +8,  x^2 + \tfrac{7}{9} y^2 + \tfrac{1}{3} y}$ e $f= x +4y$}
% ANNA->JOHN:
% ho gia' calcolato la continuazione di questo esempio nella sez
% modular.  Non credo di riuscire a farlo per questo input

%%%%%%%%%%%%%%%%%%%%%%%%%%%%%%%%%%%%%
%\section{Algorithms for Computing Minimal Polynomials}
%\label{MinPoly:sec:algms}

%Let $K$ be a field, let $P = K[x_1,\dots, x_n]$, let $I$ be a
%zero-dimensional ideal, and let $f \in P$.
%\begin{remark}\label{MinPoly:rem:NFsigma}
For the basic properties of Gr\"obner bases
we refer to \citet{B3} and \citet{KR1}.
Let $\sigma$ be a term-ordering on~$\TT^n$, and let $I$ be an ideal in
the polynomial ring $P$.
It is known that $\NF_{\sigma, I}(f)$, the $\sigma$-normal
form of $f$ with respect to $I$, does not depend on which $\sigma$-Gr\"obner
basis of $I$ is used nor on which specific rewriting steps were
used to calculate it (see~Proposition~2.4.7 in~\citet{KR1}).
If $I$ is clear from the context, we write simply~$\NF_{\sigma}(f)$.
%\end{remark}

\begin{remark}[Elimination]
\label{MinPoly:rem:elimination}
A well-known method for
computing the minimal polynomial $\minpoly{f,I}(z)$ is by elimination.
One extends $P$ with a new indeterminate to produce a new polynomial
ring $R =K[x_1,\dots,x_n,\;z]$, then defines the ideal $J =
IR+\ideal{z-f}$ in $R$, and finally eliminates the indeterminates
$x_1, \dots, x_n$.

However, even in the case where $f$ is just an indeterminate,
  the algorithm which derives from this approach is usually
  impractically slow on non-trivial examples;
though the performance could be improved by using FGLM
(see Remark~\ref{MinPoly:rem:FGLM}).
\end{remark}

% Here we give a refined version of~Algorithm~5.1.1 of~\citet{KR3}.

% \begin{ABPRAlgorithm}
%  \textsc{MinPolyQuotElim}
%  \label{MinPoly:alg:MinPolyQuotElim}
%  \begin{description}[topsep=0pt,parsep=1pt]
%   \item[\textit{notation:}] Let $P =K[x_1,\dots,x_n]$
%   \item[Input] $I$, a zero-dimensional ideal in $P$, and a polynomial $f\in P$
%   \item[1] create the extended polynomial ring $R=K[x_1,\dots,x_n,\;z]$
%   \item[2] define the ideal $J=I\cdot R+\ideal{z-f}$
%   \item[3] \textbf{return} the monic generator of the principal ideal $J\cap K[z]$
%   \item[Output] $\minpoly{f,I}(z) \in K[z]$
%  \end{description}
% \end{ABPRAlgorithm}

% \medskip

Another way to compute $\minpoly{f,I}(z)$ is via multiplication
endomorphisms on $P/I$.
Let $f \in P$ then we write $\theta_{\bar f}:P/I\rightarrow P/I$
for the endomorphism ``multiplication by $\bar{f}$''.
There is a natural isomorphism between~$P/I$ and  $K[\theta_{\bar f} \mid \bar{f} \in P/I]$
which associates $\bar{f}$ with $\theta_{\bar{f}}$
(see~Proposition~4.1.2 in~\citet{KR3}).

\begin{example}[Example~\ref{MinPoly:ex:toy-Min-neChar} continued]
\label{MinPoly:ex:toy-Min-neChar.MultMat}
Let $P =\FF_{101}[x,y]$,
$I = \ideal{y^3 -xy -2y^2 +y,\;  xy^2,\; x^2 -x}$
and $f=y$.
The given generators are a DegRevLex-Gr\"obner basis,
thus a quotient basis for
$P/I$ is $QB=\{1,  y,  y^2,  x,  xy\}$, and the matrix for ``multiplication
by~$\bar{f}$'' with respect to the basis $QB$ is
${A = \left( \begin{smallmatrix}
  0 & 0 & 0 & 0 & 0 \\
  1 & 0 & -1 & 0 & 0 \\
  0 & 1 & 2 & 0 & 0 \\
  0 & 0 & 0 & 0 & 0 \\
  0 & 0 & 1 & 1 & 0
\end{smallmatrix} \right)
}$.  In \cocoa this can be computed via the call \texttt{MultiplicationMat(y, I, QB)}.
Note that the matrix is computed by calculating $\NF_I(f\cdot t)$ for all $t\in
QB$ (see Remark~\ref{MinPoly:rem:NF}).
\end{example}

The minimal polynomial of $f$ modulo $I$ is the same as the
minimal polynomial of the endomorphism $\theta_{\bar{f}}$.
Thus, if the matrix $A$ represents $\theta_{\bar{f}}$ with
respect to some $K$-basis of $P/I$, we can compute the 
minimal polynomial of $A$ (and thus of $\theta_{\bar{f}}$)
using the following algorithm which is a refined version
of~Algorithm~1.1.8 in~\citet{KR3}.

\begin{ABPRAlgorithm}
  \textsc{MinPolyQuotMat}
  \label{MinPoly:alg:MinPolyQuotMat}
  \begin{description}[topsep=0pt,parsep=1pt]
  \item[\textit{notation:}] Let $P =K[x_1,\dots,x_n]$
    with term-ordering $\sigma$
  \item[Input] $I$, a zero-dimensional ideal in $P$,
 and a polynomial $f\in P$
  \item[1] compute $GB$, a $\sigma$-Gr\"obner basis for $I$;\\
     from $GB$ compute $QB$, the corresponding monomial quotient
    basis of $P/I$\\
   (below we assume 1 is the first element in $QB$)
  \item[2] compute $A$,  the matrix representing the
    map $\theta_{\bar{f}}$ w.r.t.~$QB$
  \item[3] let $v_0 = (1\;0\;0\;\dots)^{\rm tr}$
     and $L = \{v_0\}$
  \item[4] \textit{Main Loop:} for $i=1,2,\dots,{\rm len}(QB)$ do
    \begin{description}[topsep=0pt,parsep=1pt]
    \item[4.1] let $v_i = A \cdot v_{i-1}$ \quad (hence we have $v_i=A^i\cdot v_0$)
    \item[4.2]
      is there a linear dependency $v_i = \sum_{j=0}^{i-1}c_j v_j$ with coefficients $c_j\in K$?
      \begin{description}[topsep=0pt,parsep=1pt]
      \item[4.2-yes] \textbf{return}  $\minpoly{f,I}(z)=z^i -\sum_{j=0}^{i-1}c_j z^j$
      \item[4.2-no] append~$v_i$ to~$L$
      \end{description}
    \end{description}
  \item[Output] $\minpoly{f,I}(z) \in K[z]$
  \end{description}
\end{ABPRAlgorithm}

% \begin{example}[Example~\ref{MinPoly:ex:toy-Min=Char} continued]
% Let $P =\FF_{101}[x,y]$,
% $I = \ideal{x^2-y,\;  y^2-2x-4}$,
% and $f=5x-3y$.
% Then $\minpoly{f,I}(z)$ is $ g(z)= z^4 +18z^2 +48z -23$.
% In this case $\deg(\minpoly{f,I}(z) = 4 = \dim_{\FF_{101}}(P/I)$,
% so $\minpoly{f,I}(z)$ coincides with the characteristic polynomial
% of the endomorphism $\theta_{\bar{f}}$.
% \end{example}

% \begin{remark}
% Assuming that the monomial $1$ appears as the first element of $QB$,
% we can make step \textsc{MinPolyQuotMat}-4.2 faster:
% the existence of a linear dependency among the matrices $A^i$ is equivalent to a
% linear dependency among just their \emph{first columns}.
% The reason is that the first column contains the coefficients of $\bar{f}^i\cdot1=\bar{f}^i$ w.r.t.~the basis $QB$.  \john{So rather than computing the sequence of the $B_i$, it suffices
%   to compute the sequence $v_0 = (1,0,0,\ldots,0)$ and $v_i = A \cdot v_{i-1}$ for $i >0$, and
%   look for a linear dependency among the $v_j$.}
% \end{remark}

\begin{example}[Example~\ref{MinPoly:ex:toy-Min-neChar.MultMat} continued]
\label{MinPoly:ex:toy-Min-neChar.AlgMat}
For $i=0,1,\dots$ the vector of $v_i$
comprises the coefficients of $\bar{f}^i$ with respect to the vector space basis
($QB$ in the algorithm).  Writing these vectors as columns gives
$C = \begin{smallmatrix}
  1 & 0 & 0 & 0 & 0 & 0 & \dots \\
  0 & 1 & 0 & -1 & -2 & -3 & \dots \\
  0 & 0 & 1 & 2 & 3 & 4 & \dots\\
  0 & 0 & 0 & 0 & 0 & 0 & \dots\\
  0 & 0 & 0 & 1 & 2 & 3 & \dots
\end{smallmatrix}
$.
Note that the first 4 columns are linearly independent, whereas the
first 5 admit the relation  $\bar{f}^4 =2\bar{f}^3 -\bar{f}^2$, hence the
algorithm stops at the fourth iteration returning
 $y^4 -2y^3 +y^2$.
Notice that in this instance the minimal polynomial has
degree~$<\deg_K(P/I)$, so it strictly divides the characteristic polynomial.
% is the lowest degree polynomial in $y$ belonging to $I$.
\end{example}

\medskip

There is a still more direct approach.  It comes from
considering the very
definition of minimal polynomial: we look for the
first linear dependency among the powers ${\bar{f}}^i$ in $P/I$.
Here we give a refined version
of~Algorithm~5.1.2 in~\citet{KR3}.
% \john{Questo metodo e` molto simile al metodo che calcola la successione $v_j$,
% solo che usa polinomi ``sparsi'' invece di matrici e vettori.}
% ANNA->JOHN:
% remark dopo algoritmo

\begin{ABPRAlgorithm}
 \textsc{MinPolyQuotDef}
 \label{MinPoly:alg:MinPolyQuotDef}
 \begin{description}[topsep=0pt,parsep=1pt]
  \item[\textit{notation:}] Let $P =K[x_1,\dots,x_n]$
    with term-ordering $\sigma$
  \item[Input] $I$, a zero-dimensional ideal in $P$, and a polynomial $f\in P$
 \item[1] compute $GB$, a $\sigma$-Gr\"obner basis for $I$;\\
     from $GB$ compute $QB$, the corresponding monomial quotient
    basis of $P/I$
 \item[2] let $f = \NF_{\sigma,I}(f)$
 \item[3] let $r_0=f^0\; (= 1)$ and $L = \{r_0\}$
 \item[4] \textit{Main Loop:} for $i=1,2,\dots,{\rm len}(QB)$ do
   \begin{description}[topsep=0pt,parsep=1pt]
   \item[4.1] compute $r_i = \NF_{\sigma,I}(f\cdot r_{i-1})$ \quad (hence we have $r_i = \NF_{\sigma,I}(f^i)$)
   \item[4.2]
     is there a linear dependency $r_i = \sum_{j=0}^{i-1}c_j r_j$ with coefficients $c_j\in K$?
     \begin{description}[topsep=0pt,parsep=1pt]
     \item[4.2-yes] \textbf{return}  $\minpoly{f,I}(z)=z^i -\sum_{j=0}^{i-1}c_j z^j$
     \item[4.2-no] append~$r_i$ to~$L$
     \end{description}
   \end{description}
 \item[Output] $\minpoly{f,I}(z) \in K[z]$
 \end{description}
\end{ABPRAlgorithm}

\begin{remark}[\textsc{MinPolyQuotMat} vs. \textsc{MinPolyQuotDef}]
\label{MinPoly:rem:NF}
Notice that algorithms \textsc{MinPolyQuotMat} and
\textsc{MinPolyQuotDef} essentially do the same computation:
the first using a matrix representation, and the second a
polynomial representation.
The main intrinsic difference is that the first algorithm
  computes the
  normal forms when constructing the multiplication matrix, whereas the
  second computes the normal forms in the main loop.
%   \john{Quindi questo approccio avra` (probabilmente) un vantaggio se il grado del polinomio
%   minimo e` basso}
% ANNA->JOHN:
% in realta' il calcolo della matrice di moltiplicazione e' meno
% costoso di quanto sembrerebbe
\end{remark}

\begin{example}[Example~\ref{MinPoly:ex:toy-Min-neChar.AlgMat} continued]
In this algorithm we do not compute the  multiplication matrix,
but we compute the normal forms of successive powers of $f$: hence
$r_0=1$,
$\;r_1=\NF_I(f)=y$,
$\;r_2=\NF_I(f\cdot r_1) = y^2$,
$\;r_3=\NF_I(f\cdot r_2) = xy +2y^2 -y$,
$r_4=\NF_I(f\cdot r_3) = 2xy +3y^2 -2y$,
and so on.
From these we implicitly construct the same sequence of columns $C$, and thence obtain the same relation,
$\minpoly{f,I}(y) = y^4 -2y^3 +y^2$.
\end{example}

\begin{remark}[\textsc{MinPolyQuotMat/Def} vs. \textsc{FGLM-Elimination}]
\label{MinPoly:rem:FGLM}
The performance of the naive elimination algorithm
(Remark~\ref{MinPoly:rem:elimination})
could be greatly improved using FGLM (see~\citet{FGLM}).
Given a Gr\"obner basis of a zero-dimensional ideal, we can compute an
elimination Gr\"obner basis using linear algebra on the quotient
basis.
Note that computing a Gr\"obner basis with respect to an elimination ordering
produces the desired minimal polynomial along with many other ``superfluous''
polynomials (which are necessary to complete the Gr\"obner basis, 
and which typically have ``uglier'' coefficients than the minimal polynomial).
We can modify the FGLM--elimination process to stop the Gr\"obner basis 
computation as soon as the minimal polynomial is found.  
The result would then be effectively quite similar to the
algorithms we presented above.
\end{remark}

These two algorithms, \textsc{MinPolyQuotMat} and
\textsc{MinPolyQuotDef}, are indeed quite simple and natural, but we
want to emphasize that a careful implementation is essential for making
them efficient.  The reward is performance which is dramatically
better than the naive elimination approach
(see the timings in Subsection~\ref{MinPoly:subsec:timings-finite}).

There are two crucial steps for achieving an efficient implementation.

\begin{remark}[Linear Algebra]
\label{MinPoly:rem:LinDepMill}
The common step in both algorithms is the search for a linear
dependency (in steps
\textsc{MinPolyQuotMat}-4.2 and \textsc{MinPolyQuotDef}-4.2).  We
implemented it \emph{incrementally} in CoCoALib \citep{cocoalib}
creating a C++ object called \texttt{LinDepMill}.  This
object accepts vectors one at a time, and says whether the last vector
it was given is linearly dependent on the earlier vectors; if so, then
it makes available the representation of the last vector as a linear
combination of the earlier vectors.  Internally, as new vectors are
supplied, \texttt{LinDepMill} simply builds up and stores a
row-reduced matrix, and keeps track of the corresponding linear
representations in terms of the input row-vectors.
Moreover, it is easy to implement this class in an efficient way
over a (small prime) finite field;
in CoCoALib its core is the class \texttt{LinDepFp} which uses machine integers.
In our experiments, checking for the linear dependency
now represents 1 to 3\% of the total computation time.

In contrast, the analogous check for the linear dependency over $\QQ$
is intrinsically more expensive, and represents most of the total time;
this motivates our modular approach presented in
Section~\ref{MinPoly:sec:modularapproach}.
\end{remark}

\begin{remark}[Powers and normal forms]
\label{MinPoly:rem:optimizing-NF}
When computing over $\FF_p$ the most expensive parts of the algorithms is in the
computation of the powers, and of the normal forms.

In step~\textsc{MinPolyQuotMat}-4.1
and in step~\textsc{MinPolyQuotDef}-4.1
we adopt an incremental approach, and do not compute
the powers $A^i$ and $f^i$:
this is quite a simple idea, but very important.

Having done that, and considering the computations over $\FF_p$,
the most expensive operation for \textsc{MinPolyQuotDef} is the normal
form (in step~4.1), generally taking more than 95\% of
the total time.  On the other hand, for \textsc{MinPolyQuotMat} the
most expensive operation is the multiplication (in step~4.1),
generally taking about 90--95\% of the time.
The construction of the multiplication matrix in
\textsc{MinPolyQuotMat} (in step~2)
requires several normal forms, but in our implementation,
again using an incremental approach,
this generally takes less than 10\% of the total time.
\end{remark}

% \begin{verbatim}
%   void IncrLinDep::myAppendVec(std::vector<RingElem> v)
%     {
%       const long n = len(v);
%       if (n != myVecLen) CoCoA_ERROR(ERR::IncompatDims,
%                                     "IncrLinDep::myAppendVec");
%       if (owner(v[0]) != myRing || !HasUniqueOwner(v))
%         CoCoA_ERROR(ERR::MixedRings, "IncrLinDep::myAppendVec");
%       const long ncols = len(myM);
%       myLinRelnValue.resize(ncols, zero(myRing));
%       myLinRelnValue.push_back(one(myRing));
%       for (long i=0; i < ncols; ++i)
%       {
%         const long col = myColTbl[i];
%         if (IsZero(v[col])) continue;
%         const RingElem c = v[col]/myM[i][col];
%         for (long j=0; j < len(myRowRepr[i]); ++j)
%           myLinRelnValue[j] -= c*myRowRepr[i][j];
%         for (long j=0; j < n; ++j)
%           v[j] -= c*myM[i][j];
%       }
%       for (long i=0; i < n; ++i)
%       {
%         if (IsZero(v[i])) continue;
%         myM.push_back(v);
%         myRowRepr.push_back(myLinRelnValue);
%         myColTbl.push_back(i);
%         myLinRelnValue.clear();
%         break;
%       }
%       //      return \textbf{\boldmath } myLinRelnValue;
%     }
% \end{verbatim}

\subsection{Timings: Computing Minimal Polynomials in Finite Characteristic}
\label{MinPoly:subsec:timings-finite}

% Each example is described by introducing a
% polynomial ring $P$, an ideal $I$ in $P$, and a polynomial $f$ in $P$
% which is denoted either by $\ell$ if it is linear or by $f$ if it is
% not linear.  The task is to compute $\minpoly{\ell,I}(z)$ or
% $\minpoly{f,I}(z)$, the minimal polynomial of $\bar{\ell}$ or
% $\bar{f}$ in $P/I$.

%.......................................................
In this subsection we present some timings for the computation of
minimal polynomials of elements in zero-dimensional affine
$\FF_p$-algebras.

The column \textbf{Example} gives the reference number to the examples
listed below.  The column \textbf{GB} gives the times to compute the
\texttt{DegRevLex}-Gr\"obner basis (in seconds); the columns
\textbf{Def, Mat}
%\textbf{MinPoly - Def},
% \textbf{Mat},
%and \textbf{Elim},
give the times (in seconds) of the computation of the
algorithms \ref{MinPoly:alg:MinPolyQuotDef}
and \ref{MinPoly:alg:MinPolyQuotMat},
%and \ref{MinPoly:alg:MinPolyQuotElim}
respectively.
We implemented the algorithm \textsc{MinPolyQuotMat} both in dense
and in sparse representation (currently called \texttt{MinPolyQuotDefLin}).
We give the timings for the latter,
which is up to twice as fast on the examples below.

Note that \cocoa,
whenever computing the Gr\"obner basis~$G$ of an
ideal~$I$, stores $G$ within the representation of~$I$.
Considering that $G$ is probably precomputed (for example, to check
whether~$I$ is zero-dimensional) the timings in \textbf{Def, Mat} do
not include the \textbf{GB} time.

%the symbol $\infty$ means more than $10$ minutes.
The column \textbf{deg}
gives the degree of the answer, as an indication of the complexity of
the output.

% The timings for \textbf{Mat} are for a preliminary implementation; it
% is clearly superior to elimination, but not (yet) competitive with
% \textbf{Def}.

As a comparison, we mention that the elimination algorithm in \cocoa
takes about~50 and~90 seconds on the first two examples, and more than
10 minutes on all the others (see
Remark~\ref{MinPoly:rem:elimination}). 

% \renzo{rivedere la frase seguente}
%  The implementation of
% elimination in Singular~\citet{Singular} is more efficient,
% nevertheless the algorithm \textbf{Def} in \cocoa is at least
% twice as fast on several examples.

\begin{table}[htbp]
\caption{Timings over prime finite fields}
\label{MinPoly:tab:PrimeFields} %% DOPO \caption
\begin{center}
\begin{tabular}{|ll|r|r||r|r|r|}%{cD{.}{.}{3.2}}
\hline
\multicolumn{2}{|l|}{\textbf{Example}} &
{\textbf{GB}} &
{\textbf{\boldmath{$\dim _K$}}} &
\multicolumn{3}{|c|}{\textbf{MinPoly}}
\\
 &&
 &
 &
{\textbf{deg}}&
{\textbf{Def}} &
{\textbf{Mat}}
\\
\multicolumn{2}{|c|}{} &
\multicolumn{1}{c|}{\textit{time}} &
\multicolumn{1}{c||}{} &
{}&
{\textit{time}} &
{\textit{time}}
\\
\hline
\ref{MinPoly:ex:charp-deg500}& $f_1=t$ &0.38  & 501 &501& 1.86 & 3.75 \\
                             & $f_2$ &        &     &501& 3.15 & 6.51 \\
                             & $f_3$ &        &     &500& 4.49 & 7.98 \\
\hline
\ref{MinPoly:ex:charp-split6}&  $f$     & 0.00   & 720 &720& 2.43 &12.29 \\
\hline
\ref{MinPoly:ex:1000000007-randomp}&$f$& 0.17   & 593 &590& 4.60 &10.02 \\
\hline
\ref{MinPoly:ex:sospetto}          &$f$& 0.01   & 464 &462& 0.90 & 3.10 \\
\hline
\ref{MinPoly:ex:23largeCI} & $f_1=z$&0.00    & 880 & 11& 0.00 & 0.02 \\
                           &   $f_2$&        &     &880& 1.06 &20.68 \\
\hline
\end{tabular}
\end{center}
\end{table}

\begin{example}\label{MinPoly:ex:charp-deg500}
The following is an example of a complete intersection in characteristic $101$.
Let $P=\FF_{101}[x,y,z,t]$.
Let
$g_1=   x y z t -18 z^3  +16 y^2  -28 t^2  -33 x,    $\,
$g_2=   x y z^3  -x +y,                              $\,
$g_3=   x^4  +x^2 y -21 y^3  -7 t^2  -25 z,          $\,
$g_4 =  y t^{11}  +26 x^3  +19 z                     $.
Let $I = \ideal{ g_1,\, g_2,\, g_3,\, g_4}$,
$f_1=t$,
$f_2 = 3z^4 -5y+x-t$ and
$f_3 = 3 y^4 z^2  -y^3 z t -12 z^4  -y^3  +z^2  -x$.
\end{example}

\begin{example}\label{MinPoly:ex:charp-split6}
This is an example which uses the defining ideal of the splitting algebra of a polynomial of degree 6.

We let $P = \FF_{101}[a_1, a_2, a_3, a_4, a_5, a_6]$, and for $j =1,\dots, 6$ let $s_j$ be the elementary symmetric
polynomial in the indeterminates $a_1, a_2, a_3, a_4, a_5, a_6$.
Then the ideal  ${I = \ideal{ s_1,\  s_2,\  s_3,\  s_4,\  s_5-7, \ s_6-1}}$ is the defining ideal of the
splitting algebra of the polynomial $x^6 -7x+1$. We let
%$\ell = a_1+2a_2+3a_3+4a_4+5a_5+6a_6$
$f = a_1+2a_2+3a_3+4a_4+5a_5+6a_6$
\end{example}

\begin{example}\label{MinPoly:ex:1000000007-randomp}

This is an example of a complete intersection in ``large'' characteristic $p=1000000007$.

Let $P=\FF_{p}[x,y,z,t]$, let
$g_1= x^5 y z t +z^3  -y^2  +73 t^2  -2 x, $\,
$g_2= x y z^6  -x +y,                      $\,
$g_3= 2 x^4  -x^2 y +34 y^3  -7 z t^2 ,    $\,
$g_4 =y t^4  +26 x^3  +z                   $.
Let $I = \ideal{ g_1,\, g_2,\, g_3,\, g_4}$
and $f = x^2t+5y$.
\end{example}

\begin{example}\label{MinPoly:ex:sospetto}
This is an example of a non-radical ideal in characteristic $p=101$
which is not a complete intersection.

Let $P=\FF_{p}[x,y,z]$, let
$g_1= (x^7-y-3z)^2$,
$g_2 =  xy^5 -7z^2 -2$,
$g_3 = yz^6-x-z+14$.
Then let $J_1 = \ideal{g_1, g_2, g_3}$, let $J_2 =\ideal{x,y,z}^2$,
let $I =J_1\cap J_2$, and let $f = x^2-3xy-z$.
\end{example}

\goodbreak

\begin{example}\label{MinPoly:ex:23largeCI}
This is an example in characteristic $p=23$.

Let $P=\FF_{p}[x,y,z]$,
let  $f_1=z$, $f_2= 3x-2y+5z$ and $I = \ideal{ g_1,\, g_2,\, g_3}$ where
\begin{align*}
\quad g_1 &= y^5 -7y^4 +2y^3 +11y^2 -y +5,\\
g_2 &= z^{11} +9z^{10} -9z^9 +7z^8 -8z^7 -4z^6 +9z^5 +z^4 -5z^3 +7z^2 +z +10,\\
g_3 &= x^{16} +8x^{15} -6x^{14} -8x^{13} +4x^{12} -4x^{11} +5x^{10} +8x^9 +5x^8 -4x^7+\\
    & \qquad 5x^6 +2x^5 -7x^4 +4x^3 +10x^2 +3x +8\\
\end{align*}

% Let $P=\FF_{p}[x,y,z]$, let
% $g_1= x^{16} +8x^{15} -6x^{14} -8x^{13} +4x^{12} -4x^{11} +5x^{10} +8x^9 +5x^8 -4x^7 +5x^6 +2x^5 -7x^4 +4x^3 +10x^2 +3x +8$, \
% $g_2= y^5 -7y^4 +2y^3 +11y^2 -y +5$,\
% $g_3= z^{11} +9z^{10} -9z^9 +7z^8 -8z^7 -4z^6 +9z^5 +z^4 -5z^3 +7z^2 +z +10$.

% Let $I = \ideal{ g_1,\ g_2,\ g_3}$
% and $\ell= 3x-2y+5z$.
\end{example}

%%%%%%%%%%%%%%%%%%%%%%%%%%%%%%%%%%%%%
\section{A Modular Approach for Minimal Polynomials}
\label{MinPoly:sec:modularapproach}

The topic of this section is to show how to
compute the minimal polynomial of an element of a
zero-dimensional affine $\QQ$-algebra using a modular approach.
In this section we describe the necessary tools to achieve this goal.

Modular reduction is a very well-known technique, however there is
no universal method for addressing the specific problems of
bad reduction arising in every application.
Our problem is no exception as we shall explain shortly.

Some results in this section are essentially known,
for instance Theorem~\ref{MinPoly:thm:RGBp} is similar to Theorem 1 in \citet{Wi}.
However, our idea is to stress the theoretical
importance of reduced Gr\"obner bases.
The main reason is that, given a term ordering,  the
reduced Gr\"obner basis of an ideal is unique, so that the theoretical results
depend just on the ideal, and not on the generators given.

In particular, we define the $\sigma$-denominator of an ideal
(see Definition~\ref{MinPoly:def:densigma}) and the $(p,\sigma)$-reduction
of an ideal (see Definition~\ref{def:psigmared}).
Then we describe the relevant notions of usable, good and bad primes
(see Subsection~\ref{Detection of Suitable Primes}) and, finally,
we put it all together to produce an algorithm which turns out to
perform quite well (see Subsection~\ref{Detection of Suitable Primes}).

\subsection{Reductions of Ideals modulo $p$}
\label{MinPoly:sec:IdealsModp}

The first matter to address is the following:
given an ideal $I$ in $\QQ[x_1, \dots, x_n]$ what does it mean to reduce
$I$ modulo a prime number $p$?  Since there is no homomorphism from $\QQ$
to $\mathbb F_p$, there is no immediate, universal answer to this question.
This problem has attracted a lot of attention over many years:
various approaches can be found in~\citet{Wi}, \citet{Pa}, \citet{Gr}, \citet{NY2018},
\citet{Ar}, \citet{IPS}, \citet{ABR} and~\citet{AoyNor}.
In this section we investigate the problem,
and provide a useful answer.  We let $P = \QQ[x_1, \dots, x_n]$,
and let $\sigma$ be term-ordering on the power-products in $P$.

\begin{definition}\label{MinPoly:rem:Zdelta}
Let $\delta \in \NN$ be positive, we use the symbol $\ZZ_\delta$ to
denote the \textbf{localization of $\ZZ$ by the multiplicative
set generated by $\delta$}, \ie~the subring of~$\QQ$
consisting of numbers represented by fractions of
type $\frac{a}{\delta^k}$ where $a\in \ZZ$ and $k\in \NN$.
Beware: some authors employ the notation $\ZZ_p$ to mean
the finite field of $p$ elements, or for $p$-adic numbers.
If $p$ is a prime number, we use the symbol
$\FF_p$ to denote the \textbf{finite field}~$\ZZ/p\ZZ$.

Observe that $\ZZ_\delta$ depends only on the radical $\Rad(\delta)$, \ie~the product of all primes
dividing $\delta$.  Furthermore, if $\delta_1, \delta_2 \in \NN$ are two positive integers
then $\Rad(\delta_1)$ divides $\Rad(\delta_2)$ if and only if
$\ZZ_{\delta_1}$ is a subring of $\ZZ_{\delta_2}$.

%\smallskip

% Beware: some authors ``abuse'' the notation $\ZZ_p$ to mean
% the finite field of $p$ elements, quite different from the
% localization at powers of $p$.
% \john{un guaio e` l'uso di $\ZZ_p$ per rappresentare gli interi
%   $p$-adici, simile alla localizzazione ``fuori'' dal primo $p$}
\end{definition}

We start with the following lemma (see Remark~2.1 in~\citet{NY2018})
which tells us about the denominators which can appear in a
normal form.

\begin{lemma}\label{MinPoly:lemma:samelocaliz}
Let $\delta \in \NN_+$, let $I$ be a non-zero ideal in $P$, let $G$ be
its reduced $\sigma$-Gr\"obner basis,
and let $f \in P$. Assume that
$f$ and $G$ have all coefficients in $\ZZ_\delta$.
\begin{enumerate}
\item Every intermediate step of rewriting $f$ via $G$
has all coefficients in $\ZZ_\delta$.

\item The polynomial $\NF_\sigma(f)$ has all coefficients in $\ZZ_\delta$.
\end{enumerate}
\end{lemma}

\begin{proof}
If $f=0$, the result is trivially true.  So we now assume $f \neq 0$. %, and let $\tau = \LT(f)$.
If $f$ can be reduced by $G$ then there exists $\tau \in \Supp(f)$ such that
$\tau = t\cdot \LT_\sigma(g)$ for some $g\in G$ and some power-product $t \in \TT^n$.
Let $c$ be the coefficient
of $\tau$ in $f$; by hypothesis $c \in \ZZ_\delta$.
Then the first step of rewriting gives $f_1 = f-c\cdot t\cdot g$ which has all coefficients in
$\ZZ_\delta$.  We can now repeat the same argument for rewriting~$f_1$, and so on.
%\cancel{Using \john{riferimento manca}Remark~\ref{MinPoly:rem:NFsigma}, We deduce that }
The final result,
when no further such rewriting is possible, is the normal form
of $f$, and by this same argument it has all coefficients in
$\ZZ_\delta$.  
Since $G$ is a Gr\"obner basis the normal form is reached after a finite number of
reduction steps, and the result is independent of the choice of reducer at each step.  These considerations prove both claims.
\end{proof}

%% \john{(sopra) Volendo possiamo anche dire che i primi in $\den(f)$ che non sono
%% in $\den_\sigma(I)$ compaiono in nei denominatori di $\NF(f)$ con la stessa potenza (o una potenza inferiore).}

The following example illustrates the lemma.

\begin{example}\label{MinPoly:ex:samelocaliz}
Let $P=\QQ[x,y]$, let $I =\ideal{ f_1, f_2}$ where $f_1=3x^3 -x^2 +1$, $f_2 =  x^2-y$,
and let $\sigma = \tt DegRevLex$.
The reduced $\sigma$-Gr\"obner basis of $I$ is $G = \{g_1,\, g_2,\, g_3\}$ where
$g_1 = y^2 +{\frac{1}{3}}^{\mathstrut}x -\frac{1}{9}y +\frac{1}{9},\; g_2= xy -\frac{1}{3}y +\frac{1}{3},$ and $g_3 = x^2-y$.  We let $f = y^3$ and note that  $f , g_1, g_2, g_3 \in \ZZ_3[x,y]$.
We have $\NF_\sigma(f) = -{\frac{1}{27}}^{\mathstrut}x -\frac{17}{81}y +\frac{8}{81}$,
and it is easy to check that the explicit coefficients in the equality
$$
f = \NF_\sigma(f) +
( xy +\tfrac{1}{9}x +\tfrac{1}{3}y -\tfrac{8}{27} )\, g_2  - ( y^2 +\tfrac{1}{9}y )\, g_3
$$
are the coefficients of a sequence of rewriting steps from $f$ to $\NF_\sigma(f)$.
As shown by the lemma, they all lie in $\ZZ_3$.
\end{example}

This lemma prompts us to make the following definitions.

\begin{definition}\label{MinPoly:def:densigma}
Let $P=\QQ[x_1,\dots,x_n]$.
\begin{enumerate}
\item Given a polynomial $f\in P$,
we define the \textbf{denominator of~$f$}, denoted by ${\den}(f)$,
to be $1$ if $f=0$, and otherwise
 the least common multiple of
the denominators of the coefficients of~$f$.

\item Given a non-zero ideal $I$ in $P$, with reduced
  $\sigma$-Gr\"obner basis $G_\sigma$, we define the
   \textbf{$\sigma$-denominator of~$I$} to be $\den_\sigma(I) = \lcm\{ \den(g) \,|\, g\in G_\sigma \}$.

\item The greatest common divisor of $\den_\sigma(I)$ where $\sigma$
ranges over all term-orderings is called
the \textbf{essential denominator} of $I$.
\end{enumerate}
\end{definition}

The following easy example shows that the number $\delta$ introduced
in Lemma~\ref{MinPoly:lemma:samelocaliz} depends on~$\sigma$.
% \john{spostiamo questo esempio a dopo la definizione qui sotto??  Oppure spostiamo def di ``denominator'' a subito prima del lemma sopra?? }

\begin{example}\label{MinPoly:ex:dependondelta}
Let $P = \QQ[x,y,z]$, let $I = \ideal{f}$ where $f=2x +3y+5z$.
Depending on the term-ordering chosen,  the number $\delta$ can be
$2$, $3$ or $5$.
So we have
 ${\den}(f) = 1$, 
 $\den_{\sigma}(I) = 2$ with $\sigma = \texttt{DegRevLex}$,
 and the essential denominator of~$I$ is~$1$.
\end{example}

Now we need one more definition.

\begin{definition}\label{MinPoly:def:reductionmodp}
Let $\delta$ be a positive integer, and~$p$ be a prime number
not dividing $\delta$.  We write $\pi_p$ to denote both
the canonical homomorphism $\ZZ_\delta \To \FF_p$ and its natural
``coefficientwise'' extensions to $\ZZ_\delta[x_1, \dots, x_n] \To \FF_p[x_1, \dots, x_n]$;
we call them all \textbf{reduction homomorphisms modulo~$p$}.
\end{definition}

The following theorem illustrates the importance of the
$\sigma$-denominator of an ideal.

\begin{theorem}\label{MinPoly:thm:RGBp}\textbf{(Reduction modulo $p$ of Gr\"obner Bases)}\\
Let $I$ be a non-zero ideal in $\QQ[x_1,\dots,x_n]$ with
reduced $\sigma$-Gr\"obner basis $G$.
Let $p$ be a prime number which does not divide $\den_\sigma(I)$.
\begin{enumerate}
\item  The set $\pi_p(G)$ is the reduced $\sigma$-Gr\"{o}bner basis of the ideal $\ideal{\pi_p(G)}$.

\item The set of the residue classes of the elements in $\TT^n{\setminus}\LT_\sigma(I)$ is an $\FF_p$-basis
of the quotient ring $\FF_p[x_1, \dots, x_n]/\ideal{\pi_p(G)}$.

\item For every polynomial $f\in \QQ[x_1,\ldots,x_n]$ such that  $p \notdiv \den(f)$
we have the equality
 $\pi_p(\NF_{\sigma, I}(f)) = \NF_{\sigma, \ideal{\pi_p(G)}}(\pi_p(f))$.
 \end{enumerate}
\end{theorem}

\begin{proof}
We start by proving claim~(a).
Every polynomial $g$ in $G$ is monic, so
$\pi_p(g)$ is monic and $\LT_\sigma(\pi_p(g))=\LT_\sigma(g)$.
Next we show that $\pi_p(G)$ is a reduced $\sigma$-Gr\"obner basis.
Let the elements of the Gr\"obner basis be
$G = \{g_1, \dots, g_s\}$.  Let $1\le i<j\le s$ and let
$f_0 = t_j g_i - t_i g_j$ be the $S$-polynomial of $(g_i,g_j)$. 
This $S$-polynomial rewrites to zero
via a finite number of steps of rewriting:
$f_{k+1} = f_k - c_k \cdot t_k \cdot g_{i_k}$ for $k=0,1,\ldots,r-1$.
Let $\delta = \den_{\sigma}(I)$, then
$f_0$ and every $g_i$ have all coefficients in $\ZZ_\delta$.
Lemma~\ref{MinPoly:lemma:samelocaliz} implies that each $c_k$ is in $\ZZ_\delta$ and
that all coefficients of each $f_k$ are in $\ZZ_\delta$.

We now show that the $S$-polynomial of the $p$-reduced pair $(\pi_p(g_i), \pi_p(g_j))$ rewrites
to zero via the set $\pi_p(G)$.  First we see that $\pi_p(f_0) = t_j \pi_p(g_i) - t_i \pi_p(g_j)$.
Now applying~$\pi_p$ to each rewriting step we get
$\pi_p(f_{k+1}) = \pi_p(f_k) - \pi_p(c_k) \cdot t_k \cdot \pi_p(g_{i_k})$.
If $\pi_p(c_k) {\ne} 0$, this is a rewriting step for $\pi_p(f_k)$,
otherwise ``nothing happens'' and we simply have $\pi_p(f_{k+1}) = \pi_p(f_k)$.

This shows that all the
  $S$-polynomials of $\pi_p(G)$ rewrite to zero, and hence that $\pi_p(G)$
  is a $\sigma$-Gr\"obner basis. Finally we observe that
  $\Supp(\pi_p(g_i)) \subseteq \Supp(g_i)$ for all $i =1,\dots, s$,
  hence $\pi_p(G)$ is actually the reduced $\sigma$-Gr\"{o}bner basis
  of the ideal $\ideal{\pi_p(G)}$.

\smallskip
As already observed, we have $\LT_\sigma(g_i) = \LT_\sigma(\pi_p(g_i))$
for all $i = 1, \dots, s$, hence claim~(b) follows from~(a).

For part~(c) we let $\delta = \lcm(\den(f), \den_{\sigma}(I))$.
We use the same method as in the proof of part~(a) but starting with $f_0 =f$.
Once again all rewriting steps have
coefficients in $\ZZ_\delta$, and applying $\pi_p$ to them we get either
a rewriting step for $\pi_p(f)$ or possibly a ``nothing happens'' step.
Therefore the image of the final remainder
$\pi_p(\NF(f))$ is the normal form of $\pi_p(f)$.
\end{proof}

\goodbreak
The following example illustrates some claims of the theorem.

 \begin{example}\label{MinPoly:ex:samelocalizcontinued}
We continue the discussion of Example~\ref{MinPoly:ex:samelocaliz}. We choose $p = 2$ and
 get $\ideal{ y^2 +x +y +1,\  xy +y +1, \ x^2 +y }$ as the $(p,\sigma)$-reduction of $I$.
From Theorem~\ref{MinPoly:thm:RGBp} we know that $\{y^2 +x +y +1,\  xy +y +1, \ x^2 +y\}$
is the reduced $\sigma$-Gr\"obner basis of $\ideal{\pi_p(G)}$.
\end{example}

Theorem~\ref{MinPoly:thm:RGBp}, in particular claim~(c), motivates the following definition.

\begin{definition}\label{def:psigmared}
In the context of Theorem~\ref{MinPoly:thm:RGBp}, let  the reduced $\sigma$-Gr\"obner basis $G$
be $\{g_1, g_2, \ldots, g_s\}$.
%If $p$ is a prime number which does not divide $\den_\sigma(I)$,
Then the ideal generated by
 the set $\pi_p(G)=\{\pi_p(g_1),\dots,\pi_p(g_s)\}$ in the polynomial 
 ring $\mathbb{F}_p[x_1,\dots,x_n]$ is called
 the \textbf{$(p,\sigma)$-reduction of $I$}, and will be denoted by $I_{(p,\sigma)}$.
Observe that if $I$ is zero-dimensional so is $I_{(p,\sigma)}$.
 \end{definition}

The following example shows the necessity of considering the \textit{reduced} Gr\"obner basis
 in Theorem~\ref{MinPoly:thm:RGBp}.

 \begin{example}\label{MinPoly:ex:badExample}
 Let $P = \QQ[x]$, let $a$ be the product of many primes, for instance the product of the
 first $10^{6}$  prime numbers, and let $I =\ideal{ ax,\, x^2}$. The set $S=\{ax,\, x^2\}$ is
 a Gr\"obner basis of $I$, while the set $G=\{x\}$ is the reduced Gr\"obner basis of $I$.
 Reducing~$S$ modulo $p$ where $p\,|\,a$ produces the ideal $\ideal{x^2}$, while reducing
 $G$ produces the ideal $\ideal{x}$.
 \end{example}

We conduct a more thorough investigation into $(p,\sigma)$-reductions in~\citet{ABR-Modp}.

\subsection{Detection of Suitable Primes}
\label{Detection of Suitable Primes}

For this entire subsection the ideal $I$ will be zero-dimensional,
and since in Theorem~\ref{MinPoly:thm:RGBp}.(b) we have seen that,
for a suitable prime $p$,
the set $\TT^n{\setminus} \LT_\sigma(I)$
can be mapped both to a basis of of the ring $\QQ[x_1, \dots, x_n]/I$ and also to a basis
of~$\FF_p[x_1, \dots, x_n]/I_{(p,\sigma)}$, we are motivated to provide
the following definition.

\smallskip
\begin{definition}\label{MinPoly:def:f-matrix}
  Let $P = \QQ[x_1, \ldots, x_n]$ with term-ordering $\sigma$.  Let
  $I$ be a zero-dimensio\-nal ideal in the ring $P$, with reduced
  $\sigma$-Gr\"obner basis $G$; let $\delta=\den_\sigma(I)$.  Let the
  tuple $B=(t_1, t_2, \dots, t_d) = \TT^n{\setminus}\LT_\sigma(I)$
  with elements in increasing $\sigma$-order, so necessarily $t_1=1$
  and ${d =\dim_K(P/I)}$.  We denote the natural image of $B$ in
  $\QQ[x_1, \dots, x_n]$ by $B_\QQ$ and the natural image of $B$ in
  $\mathbb F_p[x_1, \dots, x_n]$ by $B_p$.  Recall that $B_\QQ$ in
  $P/I$ is a $\QQ$-basis of monomials for $P/I$ , and by
  Theorem~\ref{MinPoly:thm:RGBp}, that $B_p$ is an $\FF_p$-basis for
  $\FF_p[x_1, \dots, x_n]/\ideal{\pi_p(G)}$ if $p$ is a prime
  number which does not divide $\delta$.

Let  $f \in \ZZ_\delta [x_1, \dots, x_n]$. % have coefficients in~$\ZZ_\delta$.
We define the \textbf{$f$-power matrix},
$M_{B_\QQ}(f,r)$, to be the $d \times (r+1)$ matrix
whose $j$-th column (for $j =1, \dots, r+1$) contains the
coordinates of $\NF_{\sigma, I}(f^{j-1})$ in the basis $B_\QQ$.
Similarly, we define the \textbf{$\pi_p(f)$-power matrix} to be $M_{B_p}(\pi_p(f),r)$ the $d \times (r+1)$ matrix
whose $j$-th column contains the
coordinates of $\NF_{\sigma, I_{(p,\sigma)}}(\pi_p(f^{j-1}))$ in the basis $B_p$.
We observe that these matrices depend on both~$\sigma$ and
the corresponding ideals.
%\john{e` anche dalla upla $B$, che deriva da $\sigma$~---~occorre dirlo??}
\end{definition}

The following proposition contains useful information about reduction of matrices.

\begin{proposition}\label{MinPoly:prop:pMultMat}
Let $f \in P$ be a polynomial and let $\delta={\den}(f)\!\cdot \den_\sigma(I)$.
\begin{enumerate}
\item For every $r$, all the entries of the matrix $M_{B_\QQ}(f,r)$ are in $\ZZ_\delta$.

\item For every $r$, we have $\pi_p(M_{B_\QQ}(f,r))=M_{B_p}(\pi_p(f),r)$ for
any prime $p \notdiv \delta$.
\end{enumerate}
\end{proposition}

\begin{proof}
  Claim~(a) follows from Lemma~\ref{MinPoly:lemma:samelocaliz} applied
  to $f^j$ for $j=0,1,\ldots,r$.  Claim~(b) follows directly from
  Theorem~\ref{MinPoly:thm:RGBp}.(c).
\end{proof}

%.......................................................
\subsubsection{Usable primes}
\label{minpoly:subsec:usableprime}

We start this subsection with an elementary result which is placed here
for the sake of completeness.

\begin{lemma}\label{MinPoly:lemma:usingGauss}
Let $f, g\in \QQ[z]$ be monic polynomials such that
$g$ divides $f$, and let $\delta \in \NN_+$.
If $f$ has coefficients in $\ZZ_\delta$ then also $g$ has coefficients in $\ZZ_\delta$.
\end{lemma}

\begin{proof}
  By hypothesis we have a factorization $f =gh$ in $\QQ[z]$ for some monic $h \in \QQ[z]$.
  Set $D_f = \den(f)$, $D_g = \den(g)$ and $D_h = \den(h)$; so each of $D_f f$, $D_g g$ and $D_h h$
  is a primitive polynomial with integer coefficients.
  By Gauss's Lemma $(D_g g) (D_h h) = D_g D_h f$ is a primitive polynomial with integer coefficients.
  Hence $D_f = \pm D_g D_h$; in particular $D_g | D_f$, and consequently $\Rad(D_g) | \Rad(D_f)$.
  Since $f \in \ZZ_{\delta}[z]$ we have $\Rad(D_f) | \Rad(\delta)$, hence also $\Rad(D_g) | \Rad(\delta)$
  which implies that $g \in \ZZ_{\delta}[z]$.
\end{proof}

We now give a proposition which tells us which primes could appear in the
denominator of a minimal polynomial.
In the following proposition we use Definition~\ref{MinPoly:def:densigma}.(c).

% \john{
% \begin{proposition}\label{MinPoly:prop:goodminpoly}
%   Let $f \in P$, and let $\delta={\den}(f) \cdot \den_\sigma(I)$. Then the minimal polynomial
%   $\minpoly{f,I}(z)$ has all coefficients in $\ZZ_\delta$.
% \end{proposition}
% }

\begin{proposition}\label{MinPoly:prop:goodminpoly}
  Let $P = \QQ[x_1, \ldots, x_n]$, let $I$ be a zero-dimensional ideal
  in $P$, and let $f$ be a polynomial in~$P$.
% Then the minimal polynomial $\minpoly{f,I}(z)$ has all coefficients in $\ZZ_\delta$;
% in particular, for any term-ordering $\sigma$ let $\delta_\sigma = \den(f) \cdot \den_\sigma(I)$
% then $\minpoly{f,I}(z)$ has all coefficients in $\ZZ_{\delta_\sigma}$.
\begin{enumerate}
\item for any term-ordering $\sigma$ let $\delta_\sigma = \den(f) \cdot \den_\sigma(I)$
then the minimal polynomial $\minpoly{f,I}(z)$ has all coefficients in $\ZZ_{\delta_\sigma}$;
\item let $\delta=\den(f) \cdot D$ where $D$ is the essential denominator of $I$,
then $\minpoly{f,I}(z)$ has all coefficients in $\ZZ_{\delta}$.
\end{enumerate}
\end{proposition}

\begin{proof}
% It suffices to show that, for every term-ordering $\sigma$,
%  if $\delta_\sigma = {\den}(f) \cdot \den_\sigma(I)$ then
% the minimal polynomial $\minpoly{f,I}(z)$ has all coefficients in $\ZZ_{\delta_\sigma}$.
We prove claim~(a).
Let $\theta_{\bar{f}}$ be the $\QQ$-endomorphism of $P/I$ given by multiplication by $\bar{f}$.
It is known that $\minpoly{f,I}(z) = \minpoly{\theta_{\bar{f}}}(z)$
(see~Remark~4.1.3.(a) in~\citet{KR3}).
Let $\charpoly{\theta_{\bar{f}}}(z)$ be the characteristic polynomial of the endomorphism $\theta_{\bar{f}}$;
by definition $\charpoly{\theta_{\bar{f}}}(z) = \det(z\, {\rm id} - \theta_{\bar{f}})$.
Next, let  $d =\dim_\QQ(P/I)$,  let $B =(1, t_2, \dots, t_d) = \TT^n\setminus \LT_\sigma(I)$, let $I_d$ be the identity
matrix of size $d$, and let $M_B(\theta_{\bar{f}})$ be the matrix which represents $\theta_{\bar{f}}$ with respect to
the basis~$B$.
Then we have $\det(z\, {\rm id} - \theta_{\bar{f}}) = \det(z\,I_d - M_B(\theta_{\bar{f}}))$.
The entries of $M_B(\theta_{\bar{f}})$ are the coefficients of the representations of $\NF_\sigma(t_i f)$
in the basis $B$ for all $t_i \in B$. They are in~$\ZZ_\delta$ by Lemma~\ref{MinPoly:lemma:samelocaliz}.
So we have proved that $\charpoly{\theta_{\bar{f}}}(z) \in \ZZ_\delta[z]$.
From the Cayley-Hamilton Theorem we deduce  that $\minpoly{\theta_{\bar{f}}}(z)$ is a
divisor of $\charpoly{\theta_{\bar{f}}}(z)$.
It follows from Lemma~\ref{MinPoly:lemma:usingGauss}
that also~$\minpoly{\theta_{\bar{f}}}(z) \in \ZZ_\delta[z]$.

Claim~(b) follows easily from claim~(a) and the definition of essential denominator.
\end{proof}
%\john{Stesso commento di prima: i primi che dividono $\den(f)$ ma non $\den_\sigma(I)$ compaiono con la stessa potenza in $\minpoly(f)$ o con una potenza inferiore.}

% \begin{example}\label{smallden}
% \red{Let $P = \QQ[x,y,z]$, let $I = \ideal{z -\tfrac{1}{3}xy -2x,\   y^2 -\tfrac{1}{3}x,\   x^2 -27y -9}$,
% and let $f = 12x^3-7y^3+z-32$.
% It turns out that the $\tt Lex$ Gr\"obner basis with $y>x>z$ is
% $\{x -3y^2,\  z -y^3 -6y^2, \ y^4 -3y -1\}$ which implies that the essential denominator is $1$. Consequently, for every $f \in P$
% with integral coefficients, the minimal polynomial
% $\minpoly{f,I}(z)$ has integral coefficients. In our case we have
% $\minpoly{f,I}(z) = z^4 -8566z^3 +24148932z^2 -24637388440z -779913177008$.}
% \end{example}

\begin{remark}\label{def:GfanDen}
  To compute the essential denominator one needs know all the possible
  reduced Gr\"obner bases of $I$, in other words one needs to compute
  the Gr\"obner Fan of $I$ (see~\citet{MR88}), but actually computing
  the fan is practicable only for ``fairly simple'' ideals.
% \john{
% In Proposition~\ref{MinPoly:prop:goodminpoly} the value of $\delta$
% depends on the term-ordering $\sigma$, but clearly the value of
% $\den(\mu_{f,I})$ does not depend on the term-ordering.  Consequently,
% we may replace the factor $\den_\sigma(I)$ with the GCD of $\den_\tau(I)$
% as $\tau$ varies over all possible term-orderings; by the theory
% of Gr\"obner fans~\citet{MR88} it suffices to consider only a finite
% number of distinct term-orderings, but actually computing the Gr\"obner
% fan is practicable only for ``fairly simple'' ideals.
% }
\end{remark}

\begin{example}\label{smallden}
Let $P = \QQ[x,y]$ and $I = \ideal{2x+3y, y^2-4}$.
There are just two possible Gr\"obner bases: $\{ x +\frac{3}{2} y,\, y^2  -4 \}$ and
$\{ y +\frac{2}{3} x,\, x^2  -9 \}$.  Hence the essential denominator is $\gcd(2,3) = 1$.

Thus we know that the minimal polynomial of any polynomial with integer coefficients
has integer coefficients.  For instance, let $f = 23x+17y$ then $\minpoly{f,I}(z) = z^2-1225$.
\end{example}

The conclusion of the proposition above motivates the following definition.

\begin{definition}\label{MinPoly:def:usableprime}
  Let $f \in P$ be a polynomial, and let $p$ be a prime number.  Then $p$ is called a
  \textbf{usable prime for $f$ with respect to $(I, \sigma)$}
  if it does not divide $\den(f)\cdot\den_\sigma(I)$.
  If~$I$ and $\sigma$ are clear from
  the context, we say simply a \textbf{usable prime}.
%  A prime which is not usable is called \textbf{ugly}.
    It follows from the definition that, for a given input  $(f, I, \sigma)$,
  there are only finitely many unusable primes, and
  it is easy to recognize and avoid them.
\end{definition}

\subsubsection{Good and Bad Primes}

In this subsection we refine the definition of usable.

\begin{definition}\label{MinPoly:def:goodprime}
  Let $p$ be a usable prime for $f$ with respect to $(I, \sigma)$;
  consequently, by Proposition~\ref{MinPoly:prop:goodminpoly}, $\pi_p(\minpoly{f,I}(z))$ is well-defined.
  We say that $p$ is a \textbf{good prime for $f$} if
  $\minpoly{\pi_p(f),I_{(p,\sigma)}}(z) = \pi_p(\minpoly{f,I}(z))$,
in other words if
the minimal polynomial of the $p$-reduction of $f$ modulo the $(p,\sigma)$-reduction of $I$
equals the $p$-reduction of the minimal polynomial of $f$ modulo $I$ over the rationals.
Otherwise, it is called \textbf{bad}.
%%if it is usable but equality does not hold.
\end{definition}

The following simple example illustrates how a prime can be bad even if it is usable.

\begin{example}
  Let $P= \QQ[x,y]$, let $I=\ideal{ x^2,\,y^2}$, and let  $f=x+y$.
  The set $\{ x^2,\, y^2 \}$ is a reduced Gr\"{o}bner basis of $I$ for every term-ordering,
   $B= (1,x,y,xy)$ is a quotient basis of $\QQ[x,y]/I$ regardless of term-ordering.
 Moreover we have $\den(f) \cdot \den_\sigma(I)= 1$ regardless of term-ordering, and thus every
prime number is usable.  Over $\QQ$ we have
$
M_B(f,3) =  \left(
\begin{smallmatrix}
1^{\mathstrut} &   0 & 0 & 0 \cr
0 &  1 &  0 & 0 \cr
0 &   1 &  0 &  0 \cr
0 &   0 &  2 &   0
\end{smallmatrix}
\right) $.
Whence we deduce that $\minpoly{f,I}(z) = z^3$.  If we change the base field to the
finite field~$\mathbb{F}_2$, we get $ M_B(\pi_2(f),3) = \left(
\begin{smallmatrix}
1^{\mathstrut} &   0 & 0 & 0 \cr
0 &  1 &  0 & 0 \cr
0 &  1 &  0 &  0 \cr
0 &  0 &  0 &   0
\end{smallmatrix}
\right)
 $
 which shows that $\minpoly{f,I} = z^2$. It is easy to see that $2$ is the
 only bad prime in this case.
 \end{example}

Next we show that there are only finitely many bad primes.

\begin{theorem}\label{MinPoly:thm:almostallgood}\textbf{(Finitely many bad primes)}\\
Let $P = \QQ[x_1, \dots, x_n]$, let $I$ be a zero-dimensional ideal in $P$,
let $\sigma$ be a term-ordering on $\TT^n$,  let $f \in P$, and
$p$ be a usable prime.
\begin{enumerate}
\item
Then $\pi_p(\minpoly{f,I}(z))$ is a multiple of
 $\minpoly{\pi_p(f), I_{(p,\sigma)}}(z)$.

\item
There are only finitely many bad primes.

\item
The prime $p$ is good if and only if
$\deg(\minpoly{\pi_p(f), I_{(p,\sigma)}}(z)) = \deg(\minpoly{f,I}(z))$.

\end{enumerate}
\end{theorem}

\begin{proof}
To simplify the presentation  we let $\mu(z) = \minpoly{f,I}(z)$
and $\mu_p(z) = \minpoly{\pi_p(f), I_{(p,\sigma)}}(z)$.
  Let $\mu(z) {=} z^r+c_{r-1}z^{r-1} +\cdots c_0$, and
  set $\delta =\den(f) \cdot \den_\sigma(I)$.
  By Proposition~\ref{MinPoly:prop:goodminpoly} we have $\mu(z) \in \ZZ_{\delta}[z]$.
  By the definition of minimal polynomial we have
  ${f^r+c_{r-1}f^{r-1} +\cdots c_0 \in I}$.  Therefore we have an equality
  $f^r+c_{r-1}f^{r-1} +\cdots c_0 = \sum h_ig_i$ for certain $h_i \in P$ where
  $\{g_1, \dots, g_s\}$ is the reduced $\sigma$-Gr\"obner basis of the ideal $I$.
  By Lemma~\ref{MinPoly:lemma:samelocaliz} we know that each $h_i \in \ZZ_\delta$.
  Since $p$ is a usable prime, it follows from Proposition~\ref{MinPoly:prop:goodminpoly}
  that we can apply~$\pi_p$ to get
  $$\pi_p(f)^r+\pi_p(c_{r-1})\,\pi_p(f)^{r-1} +\cdots + \pi_p(c_0) = \sum_{i=1}^s \pi_p(h_i) \, \pi_p(g_i)$$
  which shows that $\pi_p(\mu(f)) \in I_{(p,\sigma)}$, and hence that $\pi_p(\mu(z))$
  it is a multiple of~$\mu_p(z)$.
  So claim~(a) is proved.

To prove~(b) and~(c) it suffices to show that only a finite number of usable primes are such that
$\pi_p(\mu(z))$ is a non-trivial multiple of~$ \mu_p(z)$,
and we argue as follows.  Since~$r$ is the degree
of $\mu(z)$ we deduce that the matrix $M_{B_\QQ}(f, r-1)$ has rank $r$, hence there exists
an $r{\times} r$-submatrix of $M_{B_\QQ}(f, r-1)$ with non-zero determinant; moreover this determinant
lies in $\ZZ_{\delta}$, so can be written as $\frac{a}{\delta^s}$ for some non-zero $a \in \ZZ$
and some $s \in \NN$.  For any prime~$p$ not dividing $a \delta$, the matrix $\pi_p(M_{B_\QQ}(f,r-1))$  has maximal rank:
by Proposition~\ref{MinPoly:prop:pMultMat} we have  $\pi_p(M_{B_\QQ}(f,r-1))=M_{B_p}(\pi_p(f), r-1)$.  Hence for these primes
the degree of $\mu_p(z)$ is $r$, and the conclusion follows.
\end{proof}

The theorem tells us that bad primes are finite in number, and computations
confirm that they are very rare
(\eg~none in the examples described in
Subsection~\ref{MinPoly:sec:timings-QQ}).
There appears to be no reasonable guaranteed way to detect bad primes,
but the following corollary tells how to easily detect \textit{relatively} bad
primes.  Using this corollary we can still be misled if there are several
bad primes whose modular minimal polynomials all have the same degree (theoretically this is possible, but almost never happens in practice).

%% \cancel{However, in order to give guaranteed results,
%% the following corollary tells how to detect \textit{relatively} bad
%% primes.  Notice that we cannot detect bad primes if they \textit{all}
%% lead to minimal polynomials of the same degree: this is possible, but
%% extremely unlikely.}

% We do not have an absolute means of detecting bad primes but,
% for two different usable primes $p_1$ and $p_2$ we can compute the
% minimal polynomials of $\pi_{p_i}(f)$ with respect to
% $I_{({p_i},\sigma)}$, and by comparing degrees we can sometimes detect
% that one prime is surely bad, though without being certain that the
% other is good.

\begin{corollary}\textbf{(Detecting some bad primes)}\label{MinPoly:cor:badprimes}\\
  In the context of Theorem~\ref{MinPoly:thm:almostallgood},
  let $p_1,p_2$ be two usable primes.  Let $\mu_1 = \minpoly{\pi_{p_1}(f), I_{({p_1},\sigma)}}(z)$
  and $\mu_2 = \minpoly{\pi_{p_2}(f), I_{({p_2},\sigma)}}(z)$ be the minimal polynomials
  of the respective modular reductions.
  \begin{enumerate}
 \item  If $\deg(\mu_1) < \deg(\mu_2)$ then~$p_1$ is a bad prime.

\item  If $\deg(\mu_1) = \dim_K(P/I)$ then $p_1$ is a good prime.
\end{enumerate}
\end{corollary}
\begin{proof}
 Claim~(a) follows from parts~(a) and~(c) of Theorem~\ref{MinPoly:thm:almostallgood}.
 Claim~(b) follows from Theorem~\ref{MinPoly:thm:almostallgood}.(c) since $\dim_K(P/I)$ is an upper bound
 for the degrees of the minimal polynomials.
\end{proof}

\subsection{The Algorithm}\label{The Algorithm}
Using  the results described so far we get the following algorithm.
We emphasise a particular aspect of our implementation, the choice of $30$-bit primes
(on a $64$-bit platform):
this choice lets us use fast machine integer arithmetic while keeping the
number of iterations of the \textit{Main Loop} close to minimal.

\begin{ABPRAlgorithm}
\label{MinPoly:alg:modular}
\textsc{MinPolyQuotModular}
\begin{description}[topsep=0pt,parsep=1pt]
  \item[\textit{notation:}] $P =\QQ[x_1,\dots,x_n]$ with term-ordering $\sigma$
\item[Input]
%$P = \QQ[x_1, \dots, x_n]$, $\sigma$ a term-ordering on $\TT^n$,
$I$, a zero-dimensional ideal in $P$, and a polynomial $f \in P$
  \item[1] compute the reduced $\sigma$-Gr\"obner basis of~$I$
  \item[2] choose a usable prime $p$~---~see Definition~\ref{MinPoly:def:usableprime}.
  \item[3] compute $f_p= \pi_p(f)$ and the ideal $I_{(p,\sigma)}$.
  \item[4] compute $\mu_p= \minpoly{f_p,I_{(p,\sigma)}} \in \FF_p[z]$, the minimal polynomial of $f_p$.
  \item[5] let $\mu\crt = {\mu_p}$ \; and \; $p\crt = p$.
  \item[6] \textit{Main Loop:}
  \begin{description}[topsep=0pt,parsep=0pt]
  \item[6.1] choose a new usable prime $p$.
  \item[6.2] compute the minimal polynomial $\mu_p \in \FF_p[z]$.
  \item[6.3] if $\deg(\mu\crt) \neq \deg(\mu_p)$ then
    \begin{description}[topsep=0pt,parsep=0pt]
    \item[6.3.1] if $\deg(\mu\crt) < \deg(\mu_p)$ then
      let $\mu\crt = {\mu_p}$ \; and \; $p\crt = p$.
    \item[6.3.2] continue with next iteration of \textit{Main Loop}
    \end{description}
  \item[6.4] let $\tilde{p}\crt = p \cdot p\crt$, and let $\tilde{\mu}\crt$
    be the polynomial whose coefficients are obtained by the Chinese
    Remainder Theorem from the coefficients of~$\mu\crt$ and~$\mu_p$.
  \item[6.5] compute the polynomial
  $\mu\calc \in \QQ[z]$ whose coefficients are
  obtained as the fault-tolerant rational reconstructions of the
  coefficients of $\tilde{\mu}\crt$ modulo $\tilde{p}\crt$.
\item[6.6] were all coefficients  ``reliably'' reconstructed?
    \begin{description}[parsep=-1pt]
    \item[6.6-yes] if $\mu\calc(f)\in I$ then
      \textbf{return}  $\mu\calc$
    \item[6.6-no] \  let $\mu\crt = \tilde{\mu}\crt$ and $p\crt = \tilde{p}\crt$.
    \end{description}
\item[6.7] Continue with next iteration of \textit{Main Loop}.

\end{description}
\item[Output] $\mu\calc \in \QQ[z]$, the \textit{certified} minimal polynomial $\minpoly{f,I}$.
\end{description}
\end{ABPRAlgorithm}

% \john{\textbf{DOMANDA:} perche' controlliamo se $\mu\calc \neq 0$ in step
%   6.6-yes?  Tutti i polinomi $\mu_p$ sono monici e dello stesso grado, quindi
% la ricostruzione sara` monica, cioe` non e` nulla}
% ANNA->JOHN probabilmente era una vecchia versione del codice

\begin{proof}
The correctness and termination of this algorithm follow from
Theorem~\ref{MinPoly:thm:almostallgood}, and
Corollary~\ref{MinPoly:cor:badprimes}.
In particular, note that under our hypothesis, all bad primes
  give polynomials whose degree is \textit{too low}.  This means that
  if we have $\mu\calc(f)\in I$ (checked in step 6.6-yes) then $\mu\calc$ is
  indeed the minimal polynomial, and not a non-trivial multiple.
\end{proof}

\begin{remark}[Certified answer] \label{MinPoly:rem:Certified answer}
  When execution enters Step~6.6-yes, the value of $\mu\calc$ is
  highly likely to be correct (and will surely be so when $p\crt$
  is large enough).  Since there is nevertheless a small chance of
  the answer being wrong, either because of a wrong rational
    reconstruction, or because of a sequence of bad primes with
    compatible answers, 
we verify it by explicitly checking that $\mu\calc(f) \in I$.

Some authors, when using modular methods in this area, give algorithms where the
answer is correct ``with high probability''.  However, we want to
emphasise that our algorithm guarantees that the answer is correct.
\end{remark}

\begin{example}[Example \ref{MinPoly:ex:toy-QBigCoeff} continued]
\label{MinPoly:ex:toy-QBigCoeff-contd}
Let $P =\QQ[x,y]$,
$I = \ideal{x^2-\frac17 y^2-5,\;  y^2+4x-\frac72}$,
and $f = 3x-2y$.
Two computations modulo $p_1 = 1073741831$, and $p_2= 1073741833$ give
$\mu_{p_1} = z^4 -460175067z^3 +525914233z^2 -306782542z -191739221$
and
$\mu_{p_2} = z^4 -153391687z^3 -131478725z^2 -460174233z
-421826757$.
Notice that they have the same degree.
Their CRT combination, modulo $\tilde{p}\crt = p_1p_2$,
is
$\tilde{\mu}\crt = $
\\
{\scriptsize $z^4 -164703074540959457z^3 +352935159730627282z^2
-494109223622877543z +123527305905719987
$}
\\
(in \cocoa this is computed by \texttt{CRTPoly(mp1,p1, mp2,p2)}).
Its rational reconstruction
is $\mu\calc = z^4  +\frac{24}{7} z^3  -\frac{6527}{49} z^2
+\frac{5868}{7} z +\frac{10967}{28}$
(in \cocoa
\texttt{RatReconstructPoly($\mu\crt, p\crt$)}).
Then we verify that $\mu\calc(f) \in I$, and we may conclude that
$\mu\calc$ is indeed $\minpoly{f,I}$.
Note that the rational reconstruction algorithm requires a modulus
a bit larger than might seem necessary so that the reconstructed
values are ``reliable'' (\ie~likely correct).
%% \john{diciamo che il criterio per ``reliable'' esige un modulo piu` grande
%% del minimo assoluto che consentirebbe di ricostruire i coeff??  Altrimenti uno potrebbe pensare che bastava un solo primo per poter ricostruire...}
\end{example}

\begin{remark}[Verification]
Termination of the \textit{Main Loop} in Algorithm~\ref{MinPoly:alg:modular}
depends on the test $\mu\calc(g)\in
I$ in step \textsc{MinPolyQuotModular}-6.6-yes; however evaluating $\mu\calc(g)$
modulo $I$ is typically computationally expensive compared to the cost of a
single iteration.  For this reason, in step \textsc{MinPolyQuotModular}-6.5 we use the
fault-tolerant rational reconstruction implemented in \cocoa
(see~\citet{Abb2015}) which gives also an indication whether the
reconstructed rational is ``reliable'' (\ie~ heuristically probably
correct).  This is a computationally cheap criterion which surely indicates
``reliable'' almost as soon as $\tilde{p}\crt$ becomes large enough to
allow correct reconstruction, while also almost certainly indicating ``not
reliable'' before then.

Once a good prime has been picked, $\mu\crt$ will have the correct degree, and thereafter the degree check in step
\textsc{MinPolyQuotModular}-6.3 ensures that only results from good primes are
used; in this situation our fault-tolerant reconstruction is
equivalent to Monagan's MQRR~\citep{M04}.

Even though it has happened only extremely rarely, we have encountered ``reliable''
reconstructions which did not return the exact answer, so by default \cocoa
truly verifies that $\mu\calc(g)$ is in~$I$.
Even if performed just once this operation may be quite costly,
so we also offer a partial verification: calling \texttt{MinPoly(f,I,z,N)}
verifies for $N$ random primes $p$
that $\pi_p(\mu\calc(\pi_p(g))$ is in~$I_{(p,\sigma)}$
(in the timings table below, we check for 3 primes).
\end{remark}

\begin{remark}[No Gr\"obner basis]
\label{MinPoly:rem:NoGBasis}
A disadvantage of Algorithm~\ref{MinPoly:alg:modular} is that it needs
a Gr\"obner
basis over~$\QQ$, requiring a potentially costly computation.  We can make
a faster heuristic variant of the algorithm by working directly with the
given generators for $I$.  Let $G'$ be the set of given generators.  We
shall skip all unusable primes which divide $\den(G')$.

In steps \textsc{MinPolyQuotModular}-3 and~4 we use the ideal $\ideal{\pi_p(G')}$
instead of $I_{(p,\sigma)}$.  In the \textit{Main Loop} we skip step
\textsc{MinPolyQuotModular}-6.3, since there are no guarantees on the
degrees of bad $\mu_p$.
For instance, in Example~\ref{MinPoly:ex:badExample},
the minimal polynomial of $x$ modulo $I=\ideal{ax, x^2}$
%\john{or $I=\ideal{x^2+ax, x^2}$???  Nel mio ideale tutti gli LC sono 1.}
modulo all ``small'' primes, is $z^2$, instead of~$z$.

Thus, we keep all the $\mu_p$ but when using the chinese remainder
theorem to combine, we take only those polynomials having the same degree
as the current $\mu_p$.
%% \john{perche' non controlliamo anche gli LTI modulo $p$?  Se i LTI sono diversi, allora almeno uno dei primi e` ``bad/sospetto'', no?  Forse non importa se i LTI sono diversi a patto che i polinimi minimi abbiano
%% il giusto grado?}

In step \textsc{MinPolyQuotModular}-6.6-yes we return
directly $\mu\calc$ skipping the check that $\mu\calc(g)$ is in $I$,
since we cannot verify the answer because we want to avoid computing
Gr\"obner basis.

  There are only finitely many primes giving a bad $\mu_p$.  We can
  see this by picking some term-ordering $\sigma$, and tracing through
  the steps to compute the reduced $\sigma$-Gr\"obner basis from the
  generators $G'$.  Any prime which divides a denominator or a leading
  coefficient at any point in the computation may give a bad $\mu_p$;
  to these we add the (finitely many) bad primes for that reduced
  Gr\"obner basis.  All remaining primes will give a good $\mu_p$.
\end{remark}

 In conclusion, what do we do in \cocoa?  We recall that \cocoa,
  whenever computing the Gr\"obner basis~$G$ of an ideal~$I$, stores
  $G$ within the representation of~$I$.  This means that if $G$ has
  already been computed (\texttt{HasGBasis(I)} gives true), then we
  can use it, and the answer of \texttt{MinPoly(f,I,z)}, is fully
  guaranteed.  This happens in all the functions described in
  Section~\ref{MinPoly:sec:uses}.

In case of a direct call with no precomputed Gr\"obner basis
the implementation of \texttt{MinPoly(f,I,z)} follows
Remark~\ref{MinPoly:rem:NoGBasis} (giving a warning, if the user sets
high verbosity with \texttt{SetVerbosityLevel(80)}).
A partial verification is performed by 
\texttt{MinPoly(f,\!I,z,N)} which verifies over $N$ more primes.

In practice, we have observed no significant advantage in skipping
the computation of the Gr\"obner basis, but we keep this code for
flexibility and further investigations.

% %%%%%%%%%%%%%%%%%%%%%%%%%%%%%
% \section{Timings}
% \label{MinPoly:sec:timings}
% %%%%%%%%%%%%%%%%%%%%%%%%%%%%%

%%%%%%%%%%%%%%%%%%%%%%%%%%%%%
\subsection{Timings: Computing Minimal Polynomials over $\QQ$}
\label{MinPoly:sec:timings-QQ}
%%%%%%%%%%%%%%%%%%%%%%%%%%%%%

In this subsection we present some timings for the computation of
minimal polynomials of elements in zero-dimensional affine
$\QQ$-algebras.

The column \textbf{Example} gives the reference number to the examples
listed below.  The column \textbf{GB} gives the times to compute the
\texttt{DegRevLex}-Gr\"obner basis.  Under the heading \textbf{MinPoly}
the column
$\QQ$ gives the times of the direct
computation over $\QQ$ using Algorithm~\ref{MinPoly:alg:MinPolyQuotDef}) \textsc{MinPolyQuotDef}; under the
sub-heading \textbf{Modular} the first sub-column gives the times
of the computation using Algorithm~\ref{MinPoly:alg:modular} \textsc{MinPolyQuotModular} (which internally uses
\textsc{MinPolyQuotDef} for each modular computation) with full
verification over $\QQ$, \ie~checking that the reconstructed
polynomial actually vanishes on the input polynomial, while the second
sub-column gives the time with a heuristic verification over $\FF_p$
for 3 random primes between $10^9$ and $2\cdot 10^9$; the third sub-column
\textbf{(\#$p$)} gives the number of primes used for the
reconstruction.  In the tables, all times are in seconds.
%, and we have written $\infty$ to mean a time longer than 10~minutes.

The columns \textbf{coeff} and \textbf{deg}
give an indication of the size of the minimal polynomial:
the first expresses the maximum magnitude of the numerators
and denominators of the coefficients, and the second is the degree.

\begin{remark}
\cocoa also offers \texttt{RingTwinFloat} arithmetic, an implementation of
\textit{heur\-istically guaranteed} floating point numbers \citep{Abb2012}.
We have also tried our algorithms using this representation of~$\QQ$, but
the modular approach gave us better timings.
\end{remark}

% \renzo{rivedere il remark seguente}
% \begin{remark}
%   A comparison of the new algorithm with elimination in Singular~\citet{Singular}
%   on several of the examples below shows that \cocoa is decidedly faster.
% \end{remark}

\begin{table}[htbp]
\caption{Timings over the rationals}
\label{MinPoly:tab:rationals} %% DOPO \caption
\centering

\begin{tabular}{|c@{}c||r|r||r|r|rr| r|r|c|}
\hline
\multicolumn{2}{|c||}{\small\textbf{Example}} &
\multicolumn{1}{c|}{\small\textbf{GB}} &
\multicolumn{1}{c||}{\small\textbf{\!\!\boldmath{$\dim _K$}}} &
\multicolumn{7}{c|}{\textbf{MinPoly}}
%\multicolumn{2}{c}{\textbf{MinPoly}}
\\
\multicolumn{2}{|c||}{} &
\multicolumn{1}{c|}{} &
\multicolumn{1}{c||}{} &
{\textbf{deg}} &
\multicolumn{1}{c|}{\textbf{$\QQ$}} &
\multicolumn{4}{c|}{\textbf{Modular}} &
%%  &{(\#$p$)} &
{\textbf{coeff}}
\\
\cline{5-11}
\multicolumn{2}{|c||}{} &
\multicolumn{1}{c|}{} &
\multicolumn{1}{c||}{} &
\multicolumn{1}{c|}{} &
{{}} &
\multicolumn{2}{c|}{\small\textit{$\QQ$-verif}} &
\multicolumn{1}{c|}{\small\textit{3-verif}} &
{} &
{{}}
\\
\multicolumn{2}{|c||}{} &
{\textit{}} &
{{}} &
{{}} &
{{\small\textit{}}} &
{\small\textit{tot}} &
{\small\textit{verif}} &
{\small\textit{verif}} &
{{}} &
{{}}
\\
\multicolumn{2}{|c||}{} &
{\textit{time}} &
{{}} &
{{}} &
{{\small\textit{time}}} &
{\small\textit{time}} &
{\small\textit{time}} &
{\small\textit{time}} &
{\#$p$} &
{{}}
\\
\hline
\ref{MinPoly:ex:QQ-rand}&$f$& 0.15& 117&116&$>600$ &15.25 &4.55& 0.26& 64 &$10^{389}, 10^{188}$\\
\hline
\ref{MinPoly:ex:QQ-CI1}&$x$& 0.00& 108&107& 47.86 & 0.39 &0.05 & 0.04& 12&$10^{93},10^{0}$\\
                       &$f$&     &    &108&224.06 & 1.31 &0.11 & 0.06& 25&$10^{210},10^{0}$\\
\hline
\ref{MinPoly:ex:QQ-CI2}&$f$& 0.00& 144&144&$>600$ & 3.77 &0.21 & 0.17& 38&$10^{330}, 10^{0}$\\
\hline
\ref{MinPoly:ex:QQ-split5} &$f$& 0.00& 120&120& 45.43 & 0.67 &0.23 & 0.13  & 9  & $10^{64}, 10^{0}$\\
\hline
\ref{MinPoly:ex:QQ-split6} &$f$& 0.00& 720&720&$>600$ &172.54&14.51& 7.80 & 58 & $10^{503}, 10^{0}$\\
\hline
\ref{MinPoly:ex:QQ-largeCI}&$z$& 0.00& 230&230&233.24&   0.39&0.10  &0.09& 5  & $10^{29}, 10^{4}$  \\
\hline
\ref{MinPoly:ex:twomaxhard}&$z$& 0.42& 149&149&89.32  &  11.30 & 10.42& 0.30 & 7 &$10^{33}, 10^{19}$\\
                           &$f$& & &149&$>600$  &  18.12&12.70& 0.38 & 30 &$10^{234}, 10^{19}$\\
\hline
\ref{MinPoly:ex:twomaxsimple}&$f$& 0.33&  55& 55& 5.33   &   0.67&0.24&  0.03 & 15 &$10^{108}, 10^{12}$\\
\hline
\ref{MinPoly:ex:QQ-PrimaryNotMax}&$y$& 0.00& 378&252& 510.85& 3.45 &1.39& 1.44 &3& $10^{11}, 10^{0}$\\
                                 &$f$&     &   &252&$>600$& 20.43&2.37& 1.66 & 26&  $10^{222}, 10^{0}$\\
\hline
\end{tabular}
\end{table}

\begin{example}\label{MinPoly:ex:QQ-rand}
This is an example with no particular structure.

Let $P=\QQ[x,y,z,t]$, let
$g_1= xyzt+83x^3+73y^2-85z^2-437t$,\,
$g_2= x^3zt+z-t$,\,
$g_3 =  zt^2+76x+94y^2-324z^3-255t^4$,\,
$g_4 = y^2z+625x+26t^3$,
and let $I = \ideal{ g_1,\, g_2,\, g_3,\, g_4}$.
Let $f = t^2 +5z$.
\end{example}

The following two examples use ideals which are complete intersections;
their reduced Gr\"obner bases are straightforward to compute.

\begin{example}\label{MinPoly:ex:QQ-CI1}
Let $P=\QQ[x,y,z,t]$, let
$g_1= x^4 +83x^3 +73y^2 -85z^2 -437t$,\,
$g_2= y^3-x$,\,
$g_3 =  z^3 +z -t$,\,
$g_4 = t^3-324z^2 +94y^2 +76x$.
Let $I = \ideal{ g_1,\ g_2,\ g_3,\ g_4}$ and
% $f_1=x$ and 
$f = 2x +3y -4z+12t$.
\end{example}

\begin{example}\label{MinPoly:ex:QQ-CI2}
Let $P=\QQ[x,y,z,t]$, let
$g_1= x^4 +83z^3 +73y^2 -t^2 -437t$,\,
$g_2= y^3 -z -t$,\,
$g_3 =  z^3 +x -t$,\,
$g_4 = t^4 -12z^2  +77y^2  +15x$.
Let $I = \ideal{ g_1,\, g_2,\, g_3,\, g_4}$ and
 $f= x -3y -12z +62t$.
\end{example}

\begin{example}\label{MinPoly:ex:QQ-split5}
This is an example which uses the defining ideal of the splitting algebra of
a polynomial of degree 5.

We let $P = \QQ[a_1, a_2, a_3, a_4, a_5]$, and for $j =1,\dots, 5$
let $s_j$ be the elementary symmetric  polynomial in the indeterminates
$a_1, a_2, a_3, a_4, a_5$.
Then we introduce  the ideal  $I = \ideal{ s_1,\, s_2,\, s_3,\, s_4 +1,\, s_5-2}$
which is the defining ideal of the
splitting algebra of the polynomial $x^5 -x -2$. We let
$f = a_1+2a_2+3a_3+4a_4+5a_5$.
\end{example}

\begin{example}\label{MinPoly:ex:QQ-split6}
This is an example which uses the defining ideal of the splitting algebra of a
polynomial of degree 6 (see also~\ref{MinPoly:ex:charp-split6}).

We let $P = \QQ[a_1, a_2, a_3, a_4, a_5, a_6]$, and for $j =1,\dots, 6$ let $s_j$
be the elementary symmetric
polynomial in the indeterminates $a_1, a_2, a_3, a_4, a_5, a_6$.
The ideal  $I = \ideal{ s_1,\, s_2,\, s_3,\, s_4,\,  s_5-7,\, s_6-1}$ is the defining ideal of the
splitting algebra of the polynomial $x^6 -7x+1$. We let
$f = a_1+2a_2+3a_3+4a_4+5a_5+6a_6$.
\end{example}

% \begin{example}\label{MinPoly:ex:QQ-87points}
% We take as $I$ the vanishing ideal of $87$ random points with
% integer coordinates $x,y,z$ (each with absolute value $<$ 10000).
% Then we let $\ell = 44x +3y -z$.
% \end{example}

\begin{example}\label{MinPoly:ex:QQ-largeCI}
In this example we let
$g_1 =z^7  -3 x -y       $,
$g_2 =y^5 z -5057 x^2  -2$,
$g_3 =x^6 y -x -z +14    $,
and $I = \ideal{g_1,\,g_2,\,g_3}$.
%  Finally we let $f=z$.
\end{example}

\begin{example}\label{MinPoly:ex:twomaxhard}
In this example we let
$f_1 = x^5 -y -3z$,
$f_2 = xy^5 -5057z^2 -2$,
$f_3 = yz^5 -x -z +14$
and $J_1 = \ideal{ f_1,\, f_2,\, f_3}$. Then we let
$g_1 = x^2 -y -3z$,
$g_2 = xy -5z^2 -12$,
$g_3 = z^3 -x -y +4$
and $J_2 = \ideal{ g_1, g_2, g_3}$. Then we let $I = J_1\cap J_2$.
We do not give explicit generators of $I$ since they are cumbersome.
Finally we let  $f = 7x -5y +2z$.
\end{example}

\begin{example}\label{MinPoly:ex:twomaxsimple}
%\john{NON CAPISCO!  Cosa cambia, cosa no, rispetto all'esempio precedente???}
This is a simplified version of Example~\ref{MinPoly:ex:twomaxhard}, in the sense
that $J_2$ and $f$ are the same.
Instead we let
$f'_1 = x^3 -y - 3z$,
$f'_2 = xy^3 - 5057z^2 -2$,
$f'_3 = yz^4 -x -z +14$
and let $J_1= \ideal{f'_1,\, f'_2,\, f'_3}$.
\end{example}

\begin{example}\label{MinPoly:ex:QQ-PrimaryNotMax}
This is an example of a non-reduced $\QQ$-algebra.
Let $P=\QQ[x,y,z]$, and
$g_1= (z^7-z-1)^2$,\
$g_2= (x^2-yz)^3$,\
$g_3 =  x^9-x-1$;
let $I = \ideal{ g_1,\ g_2,\ g_3}$
and linear form
%s $f_1 = y$ and 
$f = 2x-5y+7z$.
\end{example}

%\john{
\begin{remark}
  From Table~\ref{MinPoly:tab:rationals} it is clear that the least common denominators of
  the coefficients of the minimal polynomials are often quite small.  This is a natural
  consequence of Prop.~\ref{MinPoly:prop:goodminpoly}.(b).
\end{remark}
%}

%%%%%%%%%%%%%%%%%%%%%%%%%%%%%%%%%%%%%
\section{Uses of  Minimal Polynomials}
\label{MinPoly:sec:uses}

In this section we let $K$ be a field of characteristic zero or a perfect field of
 characteristic $p>0$ having effective $p$-th roots, and let
$I$ be a zero-dimensional ideal in the polynomial ring $P=K[x_1,\dots,x_n]$.
We start the section by recalling a useful definition.

\begin{definition}\label{MinPoly:def:squarefree}
Let $f \in P$ be a non-zero polynomial with positive degree.
We define the \textbf{square-free part}, $\sqfree(f)$, to be the product of
all distinct irreducible factors of~$f$ (which are defined up to a constant factor).
Equivalently, $\sqfree(f)$
is a  generator of the radical of the principal ideal generated by~$f$.
If~$f$ is univariate, and the coefficient field has characteristic zero then
$\sqfree(g)$ is ${\frac{g}{\gcd(g,g')}}$ up to a constant factor;
if the characteristic is positive then we can use the algorithm
described in~Proposition~3.7.12 \citep{KR1}.
\end{definition}

In the next proposition we collect important results which will be used
throughout the entire section.

\goodbreak
\begin{proposition}\label{MinPoly:prop:genericlinear}
Let $K$ be a perfect field, let $P = K[x_1, \dots, x_n]$,
let $I$ be a zero-dimensional ideal in $P$, and let $R = P/I$.
\begin{enumerate}
\item If $K$ is infinite:
\begin{enumerate}
\item[($a_1$)] If $I$ is radical,
  $\deg(\minpoly{{\ell}, I}(z)) = {\dim_K(P/I)}_{\mathstrut}$ for
 the generic linear form $\ell \in P$.

\item[($a_2$)] If $I$ is maximal, $\bar{\ell}$ is a primitive element
  of the field $P/I$ where $\ell \in P$ is the generic linear form.
\end{enumerate}

\item[(b)] If $K$ is finite:
if $I$ is maximal, there exists $f \in P$
such that $\bar{f}$ is a primitive element of the field $P/I$.
\end{enumerate}
\end{proposition}

\begin{proof}
To prove claim~(a) we observe that (a$_2$) is a special case of (a$_1$), hence
we prove (a$_1$).
Since  $I$ is radical and $K$ is infinite,  the Shape Lemma
(\eg~see~Theorem~3.7.25 in~\citet{KR1}) guarantees the existence
of  a linear change of coordinates which  brings~$I$ into normal $x_n$-position,
hence after such a transformation the last indeterminate has squarefree
minimal polynomial of degree $\dim_K(P/I)$.
Equivalently, the generic linear change of coordinates yields a situation
where the minimal polynomial of the last indeterminate is squarefree
and has degree $\dim_K(P/I)$.
If $I$ is maximal then this polynomial is necessarily irreducible.
This is exactly what
Algorithm~\ref{MinPoly:alg:IsMaximal}
tries to achieve via randomization in
step \textsc{IsMaximal}-5 (the \textit{Second Loop}).

The proof of claim~(b) follows from the well-known fact that the multiplicative group
of a finite field $L=P/I$ is cyclic, so that if $a$ is a
generator of $L\!\setminus\! \{0\}$ we have $L = K[a]$
which implies  that $\deg(\minpoly{a,I}(z)) = \dim_K(P/I)$.
We then choose $f \in P$ such that $\bar{f} = a$.
\end{proof}

The following example shows that if $K$ is finite, and $I$ is radical but not maximal then
it is possible that no element $f\in P$ exists such that $\deg(\minpoly{{\ell}, I}(z)) = {\dim_K(P/I)}$.

\begin{example}\label{MinPoly:ex:nolinearfinitefield}
Let $K = \mathbb F_2$, let  $P = K[x,y]$, $I = \ideal{x^2+x,\; y^2+y}$.
Then we can write $I$ as an intersection of maximal ideals: $I = M_1\cap M_2\cap M_3 \cap M_4$ where
$M_1 = \ideal{x,\,y}$, $M_2 = \ideal{x,\, y+1}$, $M_3 = \ideal{x+1,\, y}$,
$M_4 = \ideal{x+1,\, y+1}$.
So whatever element $f \in P$ we choose, we have $\deg(\minpoly{f,I}(z)) \le 2$
while $\dim_K(P/I) = 4$.
\end{example}

%.......................................................
\subsection{IsRadical and Radical for a  Zero-Dimensional Ideal}
\label{MinPoly:subsec:radical}

The goal of this subsection is to describe algorithms for checking if $I$
is radical, and for computing the radical of $I$.
We need the following results.

\begin{proposition}\label{MinPoly:prop:generalSeidenberg}
 Let $K$ be a perfect field, let $P=K[x_1,\dots,x_n]$, let $I$ be
 a zero-dimensional ideal in $P$,  let $f_i(x_i)$ be such
 that $ I\cap K[x_i]=\ideal{ f_i(x_i)}$ for  $i=1,\dots,n$, and
let $g_i=\sqfree(f_i(x_i))$. Then we have the equality
$\sqrt{I}=I+\ideal{ g_1,\dots,g_n}$.
\end{proposition}

\begin{proof}
By~Proposition~3.7.1 in~\citet{KR1},   the polynomials $f_i(x_i)$ are non-zero.
Since the ideal $J=I+\ideal{ g_1,\dots,g_n}$
satisfies $I\subseteq J\subseteq\sqrt{I}$, we have $\sqrt{J}=\sqrt{I}$.
By Proposition 3.7.9 in~\citet{KR1} we have $\gcd(g_i,g'_i)=1$ for all $i=1,\dots, n$,
hence the conclusion follows from Seidenberg's Lemma
(see~Proposition~3.7.15 in~\citet{KR1}).
\end{proof}

Since $I\cap K[x_i] = \ideal{ \minpoly{x_i,I}(x_i)}$, the above proposition
can be rewritten as follows.

\begin{corollary}\label{MinPoly:cor:radicalconstruction}
 Let $K$ be a perfect field, let $P=K[x_1,\dots,x_n]$, let $I$ be
 a zero-dimen\-sional ideal in $P$,  and
let $g_i=\sqfree(\minpoly{x_i,I}(x_i))$.
Then we have $\sqrt{I}=I+\ideal{ g_1,\dots,g_n}$.
\end{corollary}

The following proposition shows that in some cases it is particularly easy
to show that an ideal is radical.

\begin{proposition}\label{MinPoly:prop:radical}
 Let $K$ be a perfect field, let $I$ be a zero-dimensional ideal in
 the polynomial ring $P=K[x_1,\dots,x_n]$, and let $f\in P$.
 If the polynomial $\minpoly{f,I}(z)$ is squarefree
and $\deg(\minpoly{f,I}(z))=\dim_K(P/I)$
then~$I$ is a radical ideal.
% \end{enumerate}
\end{proposition}

\begin{proof}
Consider the $K$-algebra homomorphism $\alpha_f\colon K[z]/\ideal{ \minpoly{f,I}(z)}  \to P/I$
sending $\bar{z}\!\mapsto\! \bar{f}$ which  is  injective by definition.
Since $\dim_K(K[z]/\ideal{ \minpoly{f,I}(z)})\!=\deg(\minpoly{f,I}(z))\!=\dim_K(P/I)$,
then $\alpha_f$ is also surjective and hence an isomorphism.
By assumption, the polynomial $\minpoly{f,I}(z)$ is squarefree and hence
$P/I \cong K[z]/\ideal{ \minpoly{f,I}(z)}$
is a reduced $K$-algebra which means  that $I$ is a radical ideal.
\end{proof}

The following algorithm determines whether a zero-dimensional ideal is
radical.

\goodbreak
\begin{ABPRAlgorithm}
\textsc{IsRadical0Dim}
  \label{MinPoly:alg:IsRadical0Dim}
  \begin{description}[topsep=0pt,parsep=1pt]
  \item[\textit{notation:}]
    $K$ a perfect field and $P =K[x_1,\dots,x_n]$
  \item[Input]  $I$, a zero-dimensional ideal in $P$
  \item[1] compute $d = \dim_K(P/I)$
  \item[2] \textit{Main Loop:} for $i=1,\dots,n$ do
    \begin{description}[topsep=0pt,parsep=1pt]
    \item[2.1] compute $\mu = \minpoly{x_i,I}$
    \item[2.2] if $\mu$ is not square-free then \textbf{return} \textit{false}
    \item[2.3] if $\deg(\mu)=d$ then \textbf{return} \textit{true}
    \end{description}
  \item[3] \textbf{return} \textit{true}
  \item[Output] \textit{true/false} indicating whether $I$ is radical or not.
  \end{description}
\end{ABPRAlgorithm}

\begin{proof}
Clearly, the algorithm ends after a finite number of steps and its
correctness follows from Corollary~\ref{MinPoly:cor:radicalconstruction} and
Proposition~\ref{MinPoly:prop:radical}
\end{proof}

Similarly, we have an algorithm for computing the radical of a
zero-dimen\-sional ideal.

\goodbreak
\begin{ABPRAlgorithm}
\textsc{Radical0Dim}
  \label{MinPoly:alg:Radical0Dim}

  \begin{description}[topsep=0pt,parsep=1pt]
  \item[\textit{notation:}]
    $K$ a perfect field and $P =K[x_1,\dots,x_n]$
  \item[Input] $I$, a zero-dimensional ideal in $P$
  \item[1] let $J = I$ and compute $d = \dim_K(P/J)$
  \item[2] \textit{Main Loop:} for $i=1,\dots,n$ do
    \begin{description}[topsep=0pt,parsep=1pt]
    \item[2.1] compute $\mu = \minpoly{x_i,J}$
%% WRONG!!
%% \item[2.2]  if $\deg(\mu)=d$ then return \textbf{\boldmath  $J$}
    \item[2.2] if $\mu $ is not square-free then
      \begin{description}[topsep=0pt,parsep=1pt]
      \item[2.2.1] let $\mu =\sqfree(\mu )$
      \item[2.2.2] let $J = J+\ideal{\mu (x_i)}$
      \item[2.2.3] compute $d = \dim_K(P/J)$
      \end{description}
    \item[2.3] if $\deg(\mu)=d$ then \textbf{return} $J$
    \end{description}
  \item[3] \textbf{return} $J$
  \item[Output] $J$: the radical of $I$
  \end{description}
\end{ABPRAlgorithm}

\begin{proof}
Clearly, the algorithm ends after a finite number of steps and its
correctness follows from Proposition \ref{MinPoly:prop:radical}
and Corollary~\ref{MinPoly:cor:radicalconstruction}.
\end{proof}

\begin{example}\label{MinPoly:ex:nolinearform}
  One might hope for a fast, randomized heuristic version of this
  algorithm: instead of the \textit{Main Loop} we pick a random linear
  form $\ell$, and set $\mu = \sqfree(\minpoly{\ell,J})$, and then
  update $J = J+\ideal{\mu(\ell)}$.  The example here shows that a
  single random linear form is not always sufficient; indeed, for this
  example $n$ linearly independent linear forms must be used before
  the correct result is obtained.

Let $K$ be a field, let $P = K[x_1,\ldots,x_n]$ with $n \ge 2$, and let the ideal $I = \ideal{ x_1,\ldots,x_n}^2$.
Now let $\ell \in P$ be any non-zero linear form.  Clearly $\minpoly{\ell,I}(z) = z^2$.
Hence $2=\deg(\minpoly{\ell,I}(z) )<\dim(P/I) = n+1$, and
adding $\ideal{\ell}$ to $I$ does not yield $\ideal{x_1,\ldots,x_n}$.
\end{example}

\begin{remark}[Timeout on Gr\"obner basis]
\label{MinPoly:rem:timeout}
In step \textsc{Radical0Dim}-2.2.3 we update the value of $d$.
In practice this step can be very costly.  Since the purpose of $d$
is to let the algorithm finish before having completed all iterations
of the \textit{Main Loop}, and the since update in step~2.2.3
must reduce the value of $d$, we can safely skip the update if the computation
of $\dim_K(P/J)$ takes too long:
in our implementation in CoCoALib we have set a heuristic time limit for this
step.
In the worst case the algorithm
simply performs all iterations, even though theoretically it may have
been able to stop at an earlier iteration.
The time limit we chose is
equal to half the expected time for all the remaining minimal polynomials,
based on the average time for the first ones.
\end{remark}

\begin{example}
Let $I = \ideal{x^3  +2 x^2 y -2 y z^2 ,
  5 y^4  -4 y^3 z +3 y z^3 ,
  5 z^4  +3 x y^2  -8 y z^2 }\in\QQ[x,y,z]$.
We compute the radical of~$I$ quite quickly, $\approx0.1$s, by adding
$\sqfree(\minpoly{x,I}(x))$, $\sqfree(\minpoly{y,I}(y))$, $\sqfree(\minpoly{z,I}(z))$.
However, its Gr\"obner basis is much harder to compute ($\approx10$s):\\
{\tiny $
\{
 \ x^3  -\frac{6228}{3125} x^2  -\frac{6488}{375} x y -\frac{17648}{375} y^2  +\frac{5502}{625} x z +\frac{67124}{1875} y z -\frac{1334032}{28125} z^2  -\frac{17208}{15625} x +\frac{11169272}{140625} y -\frac{120856}{3125} z,\\
x^2 y +\frac{6903}{6250} x^2  +\frac{3469}{375} x y +\frac{8224}{375} y^2  -\frac{2976}{625} x z -\frac{66881}{3750} y z +\frac{681758}{28125} z^2  +\frac{2079}{31250} x -\frac{5702572}{140625} y +\frac{12037}{625} z,\\
x y^2  -\frac{6903}{12500} x^2  -\frac{623}{250} x y -\frac{2368}{375} y^2  +\frac{1488}{625} x z +\frac{7637}{1500} y z -\frac{36997}{9375} z^2  +\frac{75909}{62500} x +\frac{294194}{46875} y -\frac{36889}{6250} z,\\
y^3  -\frac{297}{25000} x^2  -\frac{177}{500} x y -\frac{128}{125} y^2  -\frac{144}{625} x z +\frac{119109}{125000} y z +\frac{110991}{156250} z^2  +\frac{166419}{3125000} x +\frac{36036}{390625} y +\frac{2673}{312500} z,\\
x^2 z +\frac{1116}{625} x^2  +\frac{212}{25} x y -\frac{8696}{675} y^2  +\frac{1168}{375} x z +\frac{41912}{3375} y z +\frac{39034}{1875} z^2  -\frac{1512}{3125} x -\frac{339056}{9375} y -\frac{17348}{1125} z,\\
x y z -\frac{558}{625} x^2  -\frac{106}{25} x y -\frac{64}{25} y^2  +\frac{72}{125} x z +\frac{44}{25} y z -\frac{14368}{1875} z^2  +\frac{5124}{3125} x +\frac{117536}{9375} y -\frac{1158}{625} z,\\
y^2 z +\frac{54}{625} x^2  +\frac{3}{25} x y -\frac{32}{25} y^2  -\frac{36}{125} x z +\frac{162}{3125} y z +\frac{21552}{15625} z^2  -\frac{30258}{78125} x -\frac{52416}{78125} y -\frac{972}{15625} z,\\
x z^2  -\frac{27}{125} x^2  -\frac{22}{5} x y -\frac{16}{5} y^2  +\frac{18}{25} x z +\frac{247}{100} y z -\frac{2773}{375} z^2  +\frac{81}{2500} x +\frac{22832}{1875} y -\frac{66}{25} z,\\
y z^2  +\frac{27}{250} x^2  +\frac{3}{5} x y -\frac{8}{5} y^2  -\frac{9}{25} x z +\frac{81}{1250} y z +\frac{1638}{3125} z^2  -\frac{15129}{31250} x -\frac{13104}{15625} y -\frac{243}{3125} z,\\
z^3  +\frac{27}{200} x^2  +\frac{3}{4} x y -\frac{1519}{1000} y z
+\frac{819}{1250} z^2  -\frac{15129}{25000} x -\frac{3276}{3125} y
-\frac{243}{2500} z\
\}$}

This shows the advantage of using the CoCoALib timeout mechanism to
interrupt step~\textsc{Radical0Dim}-2.2.3 when it takes too long.
\end{example}

%%%%%%%%%%%%%%%%%%%%%%%%%%%%%%%%
\subsection{IsMaximal, and IsPrimary for a Zero-Dimensional Ideal}
\label{MinPoly:subsec:IsMaximal}

In this subsection, we describe methods for checking if a zero-dimensional ideal
is primary or is maximal.  To do this
we  use different strategies depending on the characteristic of the base field.
In particular, when $K$ is a finite field with $q$ elements we can use a
specific tool, namely a $K$-vector subspace of $R=P/I$,
called the \textbf{Frobenius space of $R$}.
The main property is that its dimension is exactly the number of
primary components of $I$.
For the definition and basic properties of Frobenius spaces we
refer to~Section~5.2 in~\citet{KR3}.
For convenience, we recall the definition here.

% \john{CONFUSO!  Vogliamo usare il prefisso $q$- per Frobenius?  Sopra
%   non l'abbiamo usato, ma nell definizione si`.  Perche'?  Il valore di
% $q$ non e` evidente dal campo $K$?}
% ANNA->JOHN potrebbe essere p invece di q

\begin{definition}\label{MinPoly:def:Frob}
Let $K$ be a finite field with characterstic $p$, and let
$q=p^e$ be a power of $p$.   Let $R=P/I$ be a
zero-dimensional $K$-algebra.
\begin{enumerate}
\item The map $\Phi_q:R \to R$ defined by $a \mapsto a^q$ is a $K$-linear
endomorphism of $R$ called the $q$-\textbf{Frobenius endomorphism} of $R$.

\item The fixed-point space of $R$ with respect to $\Phi_q$, namely the set
$\{ f \in R \ | \ f^q-f = 0\}$, is called the $q$-\textbf{Frobenius
  space} of $R$,
and is denoted by $\Frob_q(R)$.
\end{enumerate}
\end{definition}

\begin{remark}
Note that here we define the generalized Frobenius endomorphism $a \mapsto a^q$ instead of the classical Frobenius endomorphism, $a \mapsto a^p$.
The generalized endomorphism is just the classical
endomorphism iterated.  In this article we shall always take $q = \#K$.
\end{remark}

% \john{continuo ad essere perplesso dal prefisso $q$- perche' qui e`
%   fissato uguale alla dimensione del campo.  Forse la dimensione del campo deve essere $q^k$ per qualche intero positivo $k$?}

The following proposition describes some features  of minimal polynomials
when the zero-dimensional ideal $I$ is primary or maximal.

\begin{proposition}\label{MinPoly:prop:poweroflin}
Let $I$ be a zero-dimensional ideal in $P=K[x_1,\dots, x_n]$.
\begin{enumerate}
\item[(a)]
If $I$ is primary then for any $f \in P$ its minimal polynomial $\minpoly{f,I}(z)$
is a power of an irreducible polynomial.

\item[(b)]
If $I$ is maximal then for any $f \in P$ its minimal polynomial $\minpoly{f,I}(z)$ is irreducible.
\end{enumerate}
\end{proposition}

\begin{proof}
We use the same argument as in the proof of Proposition~\ref{MinPoly:prop:radical},
so that we get an in injective $K$-algebra homomorphism
$K[z]/\ideal{ \minpoly{f,I}(z)} \to P/I$.

Now we prove claim~(a).
If $I$ is primary then the only zero-divisors of $P/I$ are nilpotent,
hence the same property is shared by $K[z]/\ideal{ \minpoly{f,I}(z)}$
which implies that~$\minpoly{f,I}(z)$ is a power of an irreducible element.

Analogously, if $I$ is maximal then $P/I$ is a field, hence
$K[z]/\ideal{ \minpoly{f,I}(z)}$ is an integral domain which concludes the proof.
\end{proof}

We have a sort of converse of the above proposition.

\begin{proposition}\label{MinPoly:prop:isprimaryismaximal}
Let $I$ be a zero-dimensional ideal in $P=K[x_1,\dots, x_n]$,
and let $f \in P$ be
such that $\deg(\minpoly{f,I}(z)) = \dim_K(P/I)$.
\begin{enumerate}
\item[(a)] If $\minpoly{f,I}(z)$ is a power of an irreducible factor
then~$I$ is a primary ideal.

\item[(b)] If $\minpoly{f,I}(z)$ is irreducible then~$I$ is a maximal ideal.
\end{enumerate}
\end{proposition}

\begin{proof}
As in the proof of Proposition~\ref{MinPoly:prop:poweroflin} we have an injective
$K$-algebra homomorphism  $K[z]/\ideal{ \minpoly{f,I}(z)} \to P/I$.
The assumption that $\deg(\minpoly{f,I}(z)) = \dim_K(P/I)$ implies that this
endomorphism is actually an isomorphism.
Now if $K[z]/\ideal{ \minpoly{f,I}(z)}$
has only one maximal ideal, the same property is shared by $P/I$ which implies
that $I$ is a primary ideal and claim~(a) is proved.
Analogously, if $K[z]/\ideal{ \minpoly{f,I}(z)}$ is a field, then also $P/I$ is a field
which means that $I$ is a maximal ideal.
\end{proof}

%.......................................................
\subsection{IsMaximal}

Our next goal is to check whether an ideal $I$ in $P$ is maximal, and
the following algorithm provides an answer. Note that
is a true algorithm when $K$ is finite, whereas
the termination is only heuristically guaranteed
 when~$K$ is infinite.

\begin{ABPRAlgorithm}
  \label{MinPoly:alg:IsMaximal}
  \textsc{IsMaximal}
  \begin{description}[topsep=0pt,parsep=1pt]
  \item[\textit{notation:}]
    $K$ a perfect field and $P =K[x_1,\dots,x_n]$
  \item[Input] $I$, an ideal in $P$
  \item[1] if $I$ is not zero-dimensional, \textbf{return} \textit{false}
  \item[2] compute $d = \dim_K(P/I)$
  \item[3] \textit{First Loop:} for $i=1,\dots,n$ do
    \begin{description}[topsep=0pt,parsep=1pt]
    \item[3.1] compute $\mu = \minpoly{x_i,I}$
%%    \item[IsMaximalZero-2.2] factorize $g$
    \item[3.2] if $\mu$ is reducible then \textbf{return} \textit{false}
    \item[3.3] if $\deg(\mu)=d$ then \textbf{return} \textit{true}
    \end{description}
  \item[4] if $K$ is finite then
    \begin{description}[topsep=0pt,parsep=1pt]
    \item[4.1] compute $s = \dim_K(\Frob_q(P/I))$
    \item[4.2] if $s=1$ \textbf{return} \textit{true}
    else  \textbf{return}  \textit{false}
    \end{description}
  \item[5] {\small(else $K$ is infinite)}
 \textit{Second Loop:} repeat
    \begin{description}[topsep=0pt,parsep=1pt]
    \item[5.1] pick a random linear form $\ell\in P$
    \item[5.1] compute $\mu = \minpoly{\ell,I}$
    \item[5.2] if $\mu$ is reducible then \textbf{return} \textit{false}
%%    \item[IsMaximalZero-3.3] factorize $g$
    \item[5.3] if $\deg(\mu)=d$ then \textbf{return} \textit{true}
    \end{description}
  \item[Output] \textit{true/false} indicating the maximality of $I$.
  \end{description}
\end{ABPRAlgorithm}

\begin{proof}
Let us  show the correctness.
In step 3.2, if $\mu$ is reducible,
we conclude from Proposition~\ref{MinPoly:prop:poweroflin}.(b).
In step 3.3, since $\mu$ is irreducible, if $\deg(\mu)=d$ then
we conclude from Proposition~\ref{MinPoly:prop:isprimaryismaximal}.(b).
If the \textit{First Loop} completes without returning an answer,
all polynomials $\minpoly{x_i,I}(x_i)$ are irreducible and
belong to $I$, hence $I$ is radical by Seidenberg's Lemma
(see~\citet{KR1}, Proposition 3.7.15 and Corollary 3.7.16).
Now we know that $I$ is radical, we examine the two cases below.

First we consider the case when $K$ is finite.  Then the ideal~$I$ is maximal if and
only if $\dim_K(\Frob_q(P/I))=1$ (see~\citet{KR3}, Theorem 5.2.4.(b)).
Therefore, when $K$ is finite, steps 4.1 and 4.2 show that the algorithm is correct
and terminates.

Now we consider the case when $K$ is infinite.  In step 5.2 if the
minimal polynomial~$\mu$ is reducible,
Proposition~\ref{MinPoly:prop:poweroflin}.(b) tells us that $I$ is not maximal.
In step 5.3 we know that the polynomial $\mu$ is irreducible,
so if $\deg(\mu)=d$,
Proposition~\ref{MinPoly:prop:isprimaryismaximal}.(b) tells us that $I$ is maximal.
We conclude that also in this case the algorithm is correct.
Its termination follows  heuristically from Proposition~\ref{MinPoly:prop:genericlinear}.(a).
\end{proof}

Can we make this into a proper deterministic algorithm
when~$K$ is infinite?
The following remark answers this question.

\begin{remark}\label{MinPoly:rem:makedeterministic}
If $K$ is infinite, we can substitute the \textit{Second Loop} with
a check that in the special family of linear forms described in Lemma 2.1 in~\citet{Rou}
there is one whose minimal polynomial has degree~$d$.
In this way the algorithm becomes deterministic, however the coefficients
of the linear forms tend to become large and the computation more expensive.
\end{remark}

\begin{remark}\label{MinPoly:rem:genericalsoinfinite}
  Since the computation of the Frobenius space in step \textsc{IsMaximal}-4.1
  might be
  costly, one could be tempted to first try a few random linear forms
  (as in the \textit{Second Loop}).  However, our experiments show that
  computing the minimal polynomial for a random linear form has a
  computational cost very similar to that for the Frobenius space,
  while potentially furnishing less information.  In summary, there
  is no benefit from inserting such a ``heuristic step'' just before
  \textsc{IsMaximal}-4.1.
\end{remark}

%.......................................................
\subsection{IsPrimary  for a Zero-Dimensional Ideal}
The  goal of this subsection is to check whether a
zero-dimensional ideal $I$ in~$P$ is primary.
The structure of the following algorithm is very similar to the
structure of Algorithm~\ref{MinPoly:alg:IsMaximal}.
In particular, it is important to observe that
also in this case it is a true algorithm when $K$ is finite, whereas
the termination is only heuristically guaranteed
 when~$K$ is infinite.

\begin{ABPRAlgorithm}
  \label{MinPoly:alg:IsPrimary0Dim}
  \textsc{IsPrimary0Dim}

  \begin{description}[topsep=0pt,parsep=1pt]
  \item[\textit{notation:}]
    $P =K[x_1,\dots,x_n]$
  \item[Input] $I$, a zero-dimensional ideal in $P$
  \item[1] let $J = I$ and compute $d = \dim_K(P/J)$
  \item[2] \textit{First Loop:} for $i=1$ to $n$ do
    \begin{description}
    \item[2.1] compute $\mu = \minpoly{x_i,J}$
    \item[2.2] factorize $\mu$
    \item[2.3] if $\mu$ is not a power of an irreducible factor
      \textbf{return}  \textit{false}
    \item[2.4] if $\deg(\mu)=d$ then \textbf{return} \textit{true}
    \item[2.5] if $\mu$ is not square-free  then
      \begin{description}
      \item[2.5.1] let $\mu = \sqfree(\mu)$
      \item[2.5.2] let $J = J+\ideal{\mu(x_i)}$
      \item[2.5.3] compute $d = \dim_K(P/J)$
      \item[2.5.4] if $\deg(\mu)=d$ then \textbf{return}  \textit{true}
      \end{description}
    \end{description}
  \item[3] if $K$ is finite then
    \begin{description}[topsep=0pt,parsep=1pt]
    \item[3.1] compute $s = \dim_K(\Frob_q(P/I))$
    \item[3.2] if $s=1$ \textbf{return} \textit{true}
    else  \textbf{return} \textit{false}
    \end{description}
  \item[4] {\small (else $K$ is infinite)} \textit{Second Loop:} repeat
    \begin{description}[topsep=0pt,parsep=1pt]
    \item[4.1] pick a random linear form $\ell\in P$
    \item[4.2] compute $\mu = \minpoly{\ell,J}$
    \item[4.3] if $\mu$ is reducible \textbf{return} \textit{false}
    \item[4.4] if $\deg(\mu)=d$ then \textbf{return} \textit{true}
    \end{description}
  \item[Output] \textit{true/false} indicating whether $I$ is primary or not.
 \end{description}
\end{ABPRAlgorithm}

\begin{proof}
Let us show the correctness.
In the \textit{First Loop} we work with
an ideal $J$ such that
$\sqrt J = \sqrt I$, because of the change we might perform in step 2.5.
In particular $J$ is primary if and only if $I$ is primary.
Moreover, at the end of the \textit{First Loop}, $J = \sqrt I$
by Seidenberg's Lemma (see  the proof of  Algorithm~\ref{MinPoly:alg:IsMaximal}).

In step 2.3, if $\mu$ is not a power of an irreducible,
we conclude from Proposition~\ref{MinPoly:prop:poweroflin}.(a).
In step 2.4, if $\deg(\mu) = d$ and $\mu$ is a power of an irreducible
we conclude from Proposition~\ref{MinPoly:prop:isprimaryismaximal}.(a).

If the \textit{First Loop} completes without returning an answer,
we know that $J$ is radical, and now examine the two cases.

First we look at the case when $K$ is finite.
As in the case of Algorithm~\ref{MinPoly:alg:IsMaximal},
steps 3.1 and 3.2 guarantee the correctness and termination.

Now we look at the case when $K$ is infinite.
Since $J$ is radical, checking that $I$ is primary is equivalent to
checking that $J$ is maximal.
%Since the ideal $J$ is radical and from the construction it follows
%that $J = \sqrt{I}$,  to
Now, step 4 does exactly the same thing as step~5 of Algorithm~\ref{MinPoly:alg:IsMaximal}
and the proof of the correctness is the same.
Finally, the termination follows heuristically from Proposition~\ref{MinPoly:prop:genericlinear}.(a).
\end{proof}

\begin{remark}\label{MinPoly:rem:IsPrimaryTimeout}
Much as we observed in Remark~\ref{MinPoly:rem:timeout}, the computation of
$\dim_K(P/J)$ in step~\textsc{IsPrimary0Dim}-2.5.3 can safely be skipped
if it is too costly.
\end{remark}

\begin{remark}\label{MinPoly:rem:makedeterministic2}
When $K$ is infinite, to  turn this heuristically terminating
algorithm into a true algorithm
we can repeat  the observations contained
in Remark~\ref{MinPoly:rem:makedeterministic}.

\end{remark}

\begin{remark}\label{MinPoly:rem:avoidingFirstLoop}
When $K$ is finite, it would suffice to do simply steps \textsc{IsPrimary0Dim}-3.1
and \textsc{IsPrimary0Dim}-3.2 and conclude.
However, our experiments suggest that nonetheless it is often
faster to perform the \textit{First Loop}, as it is quick and frequently
determines the result.
\end{remark}

Here is an example which shows that the property of being primary
depends strongly on the base field.

\begin{example}\label{MinPoly:ex:primary over QQ} %Example 5.1.14.
Let $K$ be a field, let  $P = K[x]$, let $f(x)= x^4-10\,x^2+1$ be the minimal
polynomial of $\sqrt{2}+\sqrt{3}$,
and let $I =\ideal{f(x)}$.
Now, if $K = \QQ$, we can easily check that $f(x)$ is irreducible, hence we deduce
that $I$  is a maximal ideal.
Conversely, if $K =\FF_p$, it is known that $f(x)$ is reducible for
every prime~$p$, and hence $I$ is not a primary ideal.
\end{example}

%%%%%%%%%%%%%%%%%%%%%%%%%%%%%%%%
\subsection{Primary Decomposition  for a Zero-Dimensional Ideal}
\label{MinPoly:sec:PrimDec}

The theoretical background we shall use for computing primary decompositions of
zero-dimensional ideals in affine $K$-algebras is explained
in~Chapter~5 in~\citet{KR3}.
The main aim of this approach is to exploit our efficient algorithms for
computing minimal polynomials. Here we describe the algorithms
implemented in \cocoa.  In particular, we remark that
the algorithms for characteristic 0 (or large positive characteristic) and for finite
characteristic have the same structure except for the choice of a
partially splitting polynomial.

First we show how the partially splitting polynomial is chosen.
The function looking for a splitting polynomial has a \textit{First Loop}
over the indeterminates; if no splitting polynomial was found, it then calls the
characteristic-dependent algorithm.

In particular, if the input ideal~$I$ is primary, it returns a polynomial
$f$ such that $\minpoly{f,I}$ is a power of a single irreducible factor (together
with the token \TotalSplit).
Otherwise it returns a polynomial
$f$ such that $\minpoly{f,I}$  has at least two irreducible factors
(together with the token \PartialSplit).

The ``strange-looking'' values we return in step \textsc{PDSplittingFiniteField}-2
and step \textsc{PDSplitting}-4-yes
just emphasize that the ideal~$I$ is primary.

The three following functions reflect the implementation in \cocoa.

\begin{ABPRAlgorithm} %Algorithm 5.1.15.
 \textsc{PDSplitting}
 \label{MinPoly:alg:PDSplitting}
 \begin{description}[topsep=0pt,parsep=1pt]
  \item[\textit{notation:}]     $P =K[x_1,\dots,x_n]$
  \item[Input] $I$, a zero-dimensional ideal in $P$
   \item[1] compute $d = \dim_K(P/I)$
  \item[2] \textit{First Loop:} for $i=1,\dots,n$ do
    \begin{description}[topsep=0pt,parsep=1pt]
    \item[2.1] compute $\mu_i = \minpoly{x_i,I}$
    \item[2.2] factorize $\mu_i = \prod_j^s \mu_{ij}^{d_j}$
    \item[2.3] if $\deg(\mu_i)=d$ then \textbf{return}
      $(x_i,\; \{\mu_{ij}^{d_j} \mid j=1,\dots,s\},\; \TotalSplit)$
    \item[2.4] if $s>1$ then \textbf{return}
      $(x_i,\; \{\mu_{ij}^{d_j} \mid j=1,\dots,s\},\; \PartialSplit)$
    \end{description}
  \item[3] if $K$ is finite, \textbf{return} \textsc{PDSplittingFiniteField}$(I)$
  \item[4] \textsc{IsPrimary0Dim}$(I)$?

    \begin{description}[topsep=0pt,parsep=1pt]
    \item[4-yes] \textbf{return} $(0,\;\{z\},\; \TotalSplit)$
    \item[4-no]
      \textbf{return} \textsc{ PDSplittingInfiniteField}$(I)$
  \end{description}
  \item[Output] $(f,\; \text{factorization of }\minpoly{f,I},\; \TotalSplit/\PartialSplit)$
  \end{description}
\end{ABPRAlgorithm}

\begin{remark}
In step \textsc{PDSplitting}-4 it would be more natural to check
directly \textsc{IsMaximal}$(\sqrt{I})$, since we have already
computed all the $\mu_i$ we have practically ``for free''
$\sqrt{I} = I + \ideal{\sqfree(\mu_i) \mid i=1,\dots,n}$;
but calling \textsc{IsMaximal} entails a potentially
costly computation of a Gr\"obner basis for $\sqrt{I}$.
In our tests \textsc{IsPrimary0Dim}$(I)$ was frequently
significantly quicker.
% we already know that
% $\sqrt{I} = I + \ideal{\sqfree(\mu_i) \mid i=1,\dots,n}$,
% so we could just ask IsMaximal($\sqrt{I}$) instead of IsPrimary($I$).
% Indeed our computation of IsPrimary($I$), even though it will recompute
% all minimal polynomials, may be faster because
% it might meet an easier intermediate GB giving a condition for termination.
For instance, this is the case for Example~\ref{MinPoly:ex:QQ-PrimaryNotMax}.
\end{remark}

The following algorithm makes a good use of $\Frob_q(P/I)$, the Frobenius space of $P/I$.
Inspired by~\citet{GWW}, the detailed theoretical and computational aspects related to this concept are
described in~\citet{KR3}, Section 5.2.

\begin{ABPRAlgorithm} %Algorithm 5.1.15.
 \textsc{PDSplittingFiniteField}
 \label{MinPoly:alg:PDSplittingCharP}
 \begin{description}[topsep=0pt,parsep=1pt]
  \item[\textit{notation:}]     $P =K[x_1,\dots,x_n]$, $K$ a finite field
  \item[Input] $I$, a zero-dimensional ideal in $P$
  \item[1] compute ${\rm FrB}$ a $K$-basis of $\Frob_q(P/I)$ and let $s=\#({\rm FrB})$
  \item[2] if $s=1$ then \textbf{return}  $(0,\;\{z\},\; \TotalSplit)$
  \item[3] pick a non-constant element $f$ of the basis ${\rm FrB}$
  \item[4] compute $\mu = \minpoly{f,I}$
  \item[5] factorize $\mu = \prod \mu_j$
  \item[6] if $\deg(\mu)=s$ then \textbf{return}
    $(f,\; \{\mu_j \mid j=1,\dots,s\},\; \TotalSplit)$
  \item[7] \textbf{return}
  $(f,\; \{\mu_j \mid j=1,\dots,s\},\; \PartialSplit)$
  \item[Output] $(f,\; \text{factorization of }\minpoly{f,I},\; \TotalSplit/\PartialSplit)$
  \end{description}
\end{ABPRAlgorithm}

\begin{remark}
From Theorem~5.2.4 in \citet{KR3}
we know that for any zero-dimensional ideal $I$,
$f\in\Frob_q(P/I)$ if and only if
$\minpoly{f,I}$ factorizes into distinct linear factors with multiplicity~1.
\end{remark}

\goodbreak

\begin{ABPRAlgorithm} %Algorithm 5.1.15.
 \textsc{PDSplittingInfiniteField}
 \label{MinPoly:alg:PDSplittingChar0}
 \begin{description}[topsep=0pt,parsep=1pt]
  \item[\textit{notation:}]     $P =K[x_1,\dots,x_n]$, $K$ an infinite field
  \item[Input] $I$, a non-primary, zero-dimensional ideal in $P$
    \item[1] compute $d=\dim_K(P/I)$
    \item[2] \textit{Main Loop:} repeat:
    \begin{description}[topsep=0pt,parsep=1pt]
    \item[2.1] pick a random linear form $\ell\in P$;
    \item[2.2] compute $\mu = \minpoly{\ell,I}$
    \item[2.3] factorize $\mu = \prod_j^s \mu_j^{d_j}$
    \item[2.4] if $\deg(\mu)=d$ then \textbf{return}
      $(\ell,\; \{\mu_j^{d_j} \mid j=1,\dots,s\},\; \TotalSplit)$
    \item[2.5] if $s>1$ then \textbf{return}
      $(\ell,\; \{\mu_j^{d_j} \mid j=1,\dots,s\},\; \PartialSplit)$
    \end{description}
  \item[Output] $(\ell ,\; \text{factorization of }\minpoly{\ell,I},\; \TotalSplit/\PartialSplit)$
  \end{description}
\end{ABPRAlgorithm}

Now we are ready to see how the splittings are used to compute the
primary decomposition.

\begin{ABPRAlgorithm}
 \textsc{PrimaryDecompositionCore}
 \label{MinPoly:alg:PrimaryDecompositionCore}
 \begin{description}[topsep=0pt,parsep=1pt]
  \item[\textit{notation:}]     $P =K[x_1,\dots,x_n]$
  \item[Input] $I$, a zero-dimensional ideal in $P$
  \item[1]  let {\small $(f,\; \{\mu_j^{d_j}\mid j{=}1,\dots\!,s\},\; \TotalSplit/\PartialSplit)$} be the output of
      \textsc{PDSplitting}$(I)$
  \item[2] if $s=1$ then \textbf{return}  ($\{I\}$,\; \TotalSplit)
  \item[3] else \textbf{return}
    ($\{I{+}\ideal{\mu_j(f)^{d_j}}\mid j{=}1,\dots\!,s\}$,\; \TotalSplit/\PartialSplit)
  \item[Output] $(\{J_1,\dots,J_s\},\; \TotalSplit/\PartialSplit)$
    \quad such that $I = J_1\cap\dots\cap J_s$
 \end{description}
\end{ABPRAlgorithm}

\begin{ABPRAlgorithm}
 \textsc{PrimaryDecomposition0Dim}
 \label{MinPoly:alg:PrimaryDecomposition0Dim}
 \begin{description}[topsep=0pt,parsep=1pt]
  \item[\textit{notation:}] $P =K[x_1,\dots,x_n]$
  \item[Input] $I$, a zero-dimensional ideal in $P$
  \item[1] let $(\{J_1,\dots,J_s\},\; \TotalSplit/\PartialSplit)$\\
be the output of \textsc{PrimaryDecompositionCore}$(I)$
  \item[2] if it is {\TotalSplit},
    \textbf{return}  $\{J_1,\dots,J_s\}$
  \item[3] \textit{Main Loop:} for $i = 1,\dots,s$ do
    \begin{description}[topsep=0pt,parsep=1pt]
    \item[3.1] is $J_i$ primary?
      \begin{description}[topsep=0pt,parsep=1pt]
      \item[3.1-yes] $Dec_i = \{J_i\}$
      \item[3.1-no]  $Dec_i = \text{\textsc{PrimaryDecomposition0Dim}}(J_i)$ {\small \quad $\longleftarrow$ recursive call}
      \end{description}
    \end{description}
  \item[4]
    \textbf{return} $Dec_1 \cup \dots \cup Dec_s$
\item[Output]  the primary decomposition of $I$
 \end{description}
\end{ABPRAlgorithm}

\bigskip

The column \textbf{Example} gives the reference number to the examples listed above.
The other columns
give, respectively, the timings (in seconds) of the computation of the
algorithms~\ref{MinPoly:alg:IsRadical0Dim}, \ref{MinPoly:alg:IsMaximal},
 \ref{MinPoly:alg:IsPrimary0Dim},
 \ref{MinPoly:alg:Radical0Dim},
and~\ref{MinPoly:alg:PrimaryDecomposition0Dim},
and an indication of their answers.

% \medskip
% \red{TODO: Perch\'e manca 6.11  For the intrinsic symmetry of the
%   problem ....}

%%%%%%%%%%%%%%%%%%
\begin{table}[htbp]
\caption{Using minimal polynomials -- prime field}
\label{MinPoly:tab:using-prime} %% DOPO \caption

\centering

\begin{tabular}{|cc||r|c||r|c||r|c||r||r|r|}
\hline
\multicolumn{2}{|c||}{\textbf{Example}} &
\multicolumn{2}{c||}{\textbf{IsRadical}} &
\multicolumn{2}{c||}{\textbf{IsMaximal}} &
\multicolumn{2}{c||}{\textbf{IsPrimary}} &
\multicolumn{1}{c||}{\textbf{Radical}} &
\multicolumn{2}{c|}{\textbf{Primary Dec.}}
\\
\multicolumn{2}{|c||}{} &
\multicolumn{2}{c||}{} &
\multicolumn{2}{c||}{} &
\multicolumn{2}{c||}{} &
\multicolumn{1}{c||}{} &
\multicolumn{1}{c}{} &
\multicolumn{1}{c|}{\#Comp}
\\
\hline
\ref{MinPoly:ex:charp-deg500} && 2.51 & false& 2.50 & false& 2.52 & false & 13.70& 2.53 &  5\\
\hline
\ref{MinPoly:ex:charp-split6} && 0.02 & true & 0.00 & false& 0.00 & false & 0.02 & 1.12 & 144\\
\hline
\ref{MinPoly:ex:1000000007-randomp} && 4.37 &false & 6.63 & false & 5.95 & false &23.60 & 6.01 & 8\\
\hline
\ref{MinPoly:ex:sospetto} && 0.64 & false & 0.48 & false &  0.68 & false & 3.99 &  3.84   &  6  \\
\hline
\ref{MinPoly:ex:23largeCI}&& 0.01 & true  & 16.19 & true & 15.90 & true &0.01 &  16.10   &  1  \\
\hline
\end{tabular}
\end{table}

%%%%%%%%%%%%%%%%
\begin{table}[htbp]
\caption{Using minimal polynomials -- rationals}
\label{MinPoly:tab:using-rationals} %% DOPO \caption
\centering
\begin{tabular}{|cc||r|c||r|c||r|c||r||r|r|}
\hline
\multicolumn{2}{|c||}{\textbf{Example}} &
\multicolumn{2}{c||}{\textbf{IsRadical}} &
\multicolumn{2}{c||}{\textbf{IsMaximal}} &
\multicolumn{2}{c||}{\textbf{IsPrimary}} &
\multicolumn{1}{c||}{\textbf{Radical}} &
\multicolumn{2}{c|}{\textbf{Primary Dec.}}
\\
\multicolumn{2}{|c||}{} &
\multicolumn{2}{c||}{} &
\multicolumn{2}{c||}{} &
\multicolumn{2}{c||}{} &
\multicolumn{1}{c||}{} &
\multicolumn{1}{c}{} &
\multicolumn{1}{c|}{\#Comp}
\\
\hline
\ref{MinPoly:ex:QQ-rand} &&18.08 & false & 4.67 & false   & 4.30  & false  & 22.91 &54.33 &2\\
\hline
\ref{MinPoly:ex:QQ-CI1} && 0.46 & false   & 0.36 &  false  & 0.36 & false   & 0.49  &3.30 &2\\
\hline
\ref{MinPoly:ex:QQ-CI2} && 0.87& false  & 0.70 &false    & 0.46  & false   &1.55 & 1.11& 2\\
\hline
\ref{MinPoly:ex:QQ-split5}&&0.01 & true  & 1.38 & true   & 2.61  & true   &0.01 & 1.95 & 1\\
\hline
\ref{MinPoly:ex:QQ-split6}&& 0.03 & true  & $>600$ &   & $>600$ &  & 0.03 & $>600$ &   \\
%\hline
%\ref{MinPoly:ex:QQ-87points}&& 2.72  & true  &  1.10& false   &0.95 & false & 2.80 &  45.87& 87  \\
\hline
\ref{MinPoly:ex:QQ-largeCI} && 0.44  &true & 0.43  &true  & 0.44 &true  & 0.41 & 0.44  &  1 \\
\hline
\ref{MinPoly:ex:twomaxhard} && 11.65 &true  &12.66 &false & 12.30 & false &12.68 & 11.36 &  2\\
\hline
\ref{MinPoly:ex:twomaxsimple}&& 0.28 & true  &0.27 & false & 0.28 &false & 0.28 & 0.25  &2   \\
\hline
\ref{MinPoly:ex:QQ-PrimaryNotMax}&& 0.14 & false &0.17 & false & 1.08 &true & 1.08& 4.47  &1   \\
\hline
\end{tabular}
\end{table}

\begin{remark}
  For the example~\ref{MinPoly:ex:QQ-split6} we obtained minimal polynomials
  which are hard to factorize (like the Swinnerton-Dyer polynomials):
  they have many low degree factors modulo every prime we tried.
  The long computation times were due to the factorizer in \cocoa.
\end{remark}

% \renzo{TODO: rivedere il remark seguente}
% \begin{remark}
%   In our tests the \cocoa implementations for radical and primary decomposition are
%   generally significantly faster than those in Singular~\citet{Singular}
%   and Macaulay2~\citet{Macaulay2}.
% \end{remark}

%-------------------------------------------------------
\subsection{Comparison with Singular}

We present in Table~\ref{MinPoly:tab:SingularCoCoA} comparative timings of
our implementation with Singular \citep{Singular}; in Singular we used
the functions \texttt{radical} and \texttt{primdecGTZ}.  It is clear that
our implementation is usefully faster (and more reliable) in most cases;
an exception is the primary decomposition of example~\ref{MinPoly:ex:twomaxhard}
where the actual computation of the unverified result is fast, but
the final verification takes more than 90\% of the total time
(as shown in Table~\ref{MinPoly:tab:rationals}).  We had
hoped to include also a comparison with Macaulay2 \citep{Macaulay2},
but were unable to get timings for most of the examples.

%%%%%%%%%%%%%%%%%%
\begin{table}[htbp]
\caption{Singular and \cocoa-- time comparisons}
\label{MinPoly:tab:SingularCoCoA} %% DOPO \caption

\centering

\begin{tabular}{|c||r|r||r|r|}
\hline
\multicolumn{1}{|c||}{\textbf{Example}} &
\multicolumn{2}{c||}{\textbf{Radical}} &
\multicolumn{2}{c|}{\textbf{Primary Dec.}}
\\
\footnotesize{prime field} & Singular & \cocoa & Singular & \cocoa
\\
\hline
\ref{MinPoly:ex:charp-deg500}        &>600 & 13.70 & $^{(1)}$ 7.87 & 2.53 \\
\hline                                                   
\ref{MinPoly:ex:charp-split6}        &0.01 & 0.02  &1.00 & 1.12 \\
\hline                                                   
\ref{MinPoly:ex:1000000007-randomp}  &\tiny\texttt{\!p too large\!}&23.60  &\tiny\texttt{\!p too large\!}& 6.01\\
\hline                                                   
\ref{MinPoly:ex:sospetto}            &4.91 & 3.99  &1.62 &  3.84     \\
\hline                                                   
\ref{MinPoly:ex:23largeCI}           &0.01 &0.01   &98.15&  16.10    \\
\hline
\end{tabular}
%%%%%%%%%%%%%%%%
\begin{tabular}{|c||r|r||r|r|}
\hline
\multicolumn{1}{|c||}{\textbf{Example}} &
\multicolumn{2}{c||}{\textbf{Radical}} &
\multicolumn{2}{c|}{\textbf{Primary Dec.}}
\\
\footnotesize{rationals} & Singular & \cocoa & Singular & \cocoa
\\
\hline
\ref{MinPoly:ex:QQ-rand}         &>600   & 22.91 &151.48 &54.33\\
\hline                                                   
\ref{MinPoly:ex:QQ-CI1}          &2.17  & 0.49  &7.24  &3.30 \\
\hline                                                   
\ref{MinPoly:ex:QQ-CI2}          &40.78  &1.55   &16.41  & 1.11\\
\hline                                                   
\ref{MinPoly:ex:QQ-split5}       &0.01   &0.01   &crash  & 1.95    \\
\hline                                                   
\ref{MinPoly:ex:QQ-split6}       &0.01   & 0.03  &crash  & $>600$   \\
\hline                                                   
\ref{MinPoly:ex:QQ-largeCI}      &>600   & 0.41  &$^{(2)}$ 3.62 & 0.44     \\
\hline                                                   
\ref{MinPoly:ex:twomaxhard}      &116.29 &12.68  &3.82   & 11.36    \\
\hline                                                   
\ref{MinPoly:ex:twomaxsimple}    &0.86   & 0.28  &0.11    & 0.25     \\
\hline                                                   
\ref{MinPoly:ex:QQ-PrimaryNotMax}&93.83  & 1.08  &>600   & 4.47     \\
\hline
\end{tabular}

(1) \texttt{possible overflow}
\qquad
(2) \texttt{overflow warning}
\end{table}

%%%%%%%%%%%%%%%%%%%%%%%%%%%%%%%%%%%%%%%%%%%%%%%%%%%%%%%%%%%%%%%%%%%%%%%%%%%%%%
\section{Conclusion and future work}

We have presented both theoretical and practical aspects of our implementations
in CoCoALib for computing minimal polynomials (over $\FF_p$).
Then we presented an algorithm for computing minimal polynomials over $\QQ$ based on
a modular approach and which guarantees correctness of the result.
Finally we described several algorithms which use minimal polynomials
for various operations on zero-dimensional ideals
(\eg~testing if an ideal is radical, primary or maximal).

Our experiments have shown the potential of a good implementation, and
how this opens the way to new
applications.  For example in~\citet{AbbBigPal2018}, we use our primary
decomposition approach for factoring polynomials over algebraic field
extensions in advanced methods in the context of the {\sf SC$^2$} community: it is proving to be useful in the software CArL/SMT-RAT by~\citet{KA2018} which implements Lazard's variant of Cylindrical Algebraic Decomposition.

On the theoretical side, we investigate more deeply the consequences
of our new modular approach and apply it to general ideals
in the preprint by~\citet{ABR-Modp}.

%%%%%%%%%%%%%%%%%%%%%%%%%%%%
% Bibliography
%%%%%%%%%%%%%%%%%%%%%%%%%%%%

\bigskip\goodbreak

% Apparently required by JSC:
\bibliographystyle{elsarticle-harv} 
%\bibliographystyle{authordate} 
%\biboptions{authoryear}
\bibliography{MinPoly2019-FINAL}{}

\begin{thebibliography}{33}
\expandafter\ifx\csname natexlab\endcsname\relax\def\natexlab#1{#1}\fi
\expandafter\ifx\csname url\endcsname\relax
  \def\url#1{\texttt{#1}}\fi
\expandafter\ifx\csname urlprefix\endcsname\relax\def\urlprefix{URL }\fi

\bibitem[{Abbott(2012)}]{Abb2012}
Abbott, J., 2012. {Twin-float arithmetic}. J. Symb. Comput. 47~(5), 536--551.

\bibitem[{Abbott(2017)}]{Abb2015}
Abbott, J., 2017. {Fault-Tolerant Modular Reconstruction of Rational Numbers}.
  J.~Symb.~Comp. 80, 707--718.

\bibitem[{Abbott and Bigatti(2017)}]{AbbBig2017}
Abbott, J., Bigatti, A., 2017. {New in CoCoA-5.2.2 and CoCoALib-0.99560 for
  SC-Square}. In: Proceedings of the 2nd International Workshop on
  Satisfiability Checking and Symbolic Computation. Vol. 1974. pp. 1--6.

\bibitem[{Abbott and Bigatti(2019)}]{cocoalib}
Abbott, J., Bigatti, A., 2019. {CoCoALib:} a {C++} library for doing
  {Computations in Commutative Algebra}. Available at
  \texttt{http://cocoa.dima.unige.it/cocoalib}.

\bibitem[{Abbott et~al.(2018{\natexlab{a}})Abbott, Bigatti, and
  Palezzato}]{AbbBigPal2018}
Abbott, J., Bigatti, A., Palezzato, E., 2018{\natexlab{a}}. {New in CoCoA-5.2.4
  and CoCoALib-0.99600 for SC-Square}. In: Proceedings of the 3rd Workshop on
  Satisfiability Checking and Symbolic Computation. Vol. 2189. pp. 88--94.

\bibitem[{Abbott et~al.(2017)Abbott, Bigatti, and Robbiano}]{ABR}
Abbott, J., Bigatti, A., Robbiano, L., 2017. {Implicitization of
  Hypersurfaces}. J. Symb. Comput. 81, 20--40.

\bibitem[{Abbott et~al.(2018{\natexlab{b}})Abbott, Bigatti, and
  Robbiano}]{ABR-Modp}
Abbott, J., Bigatti, A., Robbiano, L., 2018{\natexlab{b}}. {Ideals modulo $p$}.
  \texttt{ArXiv:1801.06112}.

\bibitem[{Abbott et~al.(2019)Abbott, Bigatti, and Robbiano}]{cocoa5}
Abbott, J., Bigatti, A., Robbiano, L., 2019. {CoCoA: a system for doing
  Computations in Commutative Algebra}. Available at
  \url{http://cocoa.dima.unige.it}.

\bibitem[{Aoyama and Noro(2018)}]{AoyNor}
Aoyama, T., Noro, M., 2018. {Modular Algorithms for Computing Minimal
  Associated Primes and Radicals of Polynomial Ideals}. In: Proc.~ISSAC~2018.
  pp. 31--38.

\bibitem[{Arnold(2003)}]{Ar}
Arnold, E., 2003. {Modular algorithms for computing Gr\"obner bases}. J. Symb.
  Comput. 35, 403--419.

\bibitem[{Bostan et~al.(2003)Bostan, Salvy, and Schost}]{BSS}
Bostan, A., Salvy, B., Schost, E., 2003. {Fast Algorithms for Zero-Dimensional
  Polynomial Systems Using Duality}. AAECC 14, 239--272.

\bibitem[{Buchberger(1985)}]{B3}
Buchberger, B., 1985. {Gr\"obner Bases: An Algorithmic Method in Polynomial
  Ideal Theory}.

\bibitem[{Collins(1975)}]{collins-cad}
Collins, G., 1975. {Quantifier elimination for real closed fields by
  cylindrical algebraic decomposition}. In: Automata Theory and Formal
  Languages 2nd GI Conference. pp. 134--183.

\bibitem[{Decker et~al.(2019)Decker, Greuel, Pfister, and
  Sch{\"o}nemann}]{Singular}
Decker, W., Greuel, G.-M., Pfister, G., Sch{\"o}nemann, H., 2019. {\sc
  Singular}, \textit{{A} computer algebra system for polynomial computations}.
  Available at \texttt{http://www.singular.uni-kl.de/}.

\bibitem[{Faug\`ere et~al.(1993)Faug\`ere, Gianni, Lazard, and Mora}]{FGLM}
Faug\`ere, J., Gianni, P., Lazard, D., Mora, T., 1993. {Efficient Computation
  of Zero-dimensional Gr\"obner Bases by Change of Ordering}. J. Symb. Comput.
  16, 329--344.

\bibitem[{Gao et~al.(2008)Gao, Wan, and Wang}]{GWW}
Gao, S., Wan, D., Wang, M., 2008. {Primary Decomposition of Zero-Dimensional
  Ideals over Finite Fields}. Math. Comp. 78, 509--521.

\bibitem[{Gr\"abe(1988)}]{Gr}
Gr\"abe, H.-G., 1988. {On Lucky Primes}. J. Symb. Comput. 6, 183--208.

\bibitem[{Grayson and Stillman(2019)}]{Macaulay2}
Grayson, D., Stillman, M., 2019. Macaulay2, a software system for research in
  algebraic geometry. Available at \texttt{http://www.math.uiuc.ed/Macaulay2/}.

\bibitem[{Idrees et~al.(2011)Idrees, Pfister, and Steidel}]{IPS}
Idrees, N., Pfister, G., Steidel, S., 2011. {Parallelization of Modular
  Algorithms}. J. Symb. Comput. 46, 672--684.

\bibitem[{Kremer and \'Abrah\'am(2018)}]{KA2018}
Kremer, G., \'Abrah\'am, E., 2018. {Modular strategic SMT solving with
  SMT-RAT}.

\bibitem[{Kreuzer and Robbiano(2008)}]{KR1}
Kreuzer, M., Robbiano, L., 2008. Computational Commutative Algebra 1, 2nd
  Edition. Springer, Heidelberg.

\bibitem[{Kreuzer and Robbiano(2016)}]{KR3}
Kreuzer, M., Robbiano, L., 2016. Computational Linear and Commutative Algebra.
  Springer, Heidelberg.

\bibitem[{Lazard(1992)}]{La}
Lazard, D., 1992. {Solving zero-dimensional algebraic systems}. J. Symb.
  Comput. 13, 117--133.

\bibitem[{McCallum et~al.(2017)McCallum, Parusi\'nski, and
  Paunescu}]{lazard-CAD}
McCallum, S., Parusi\'nski, A., Paunescu, L., 2017. {Validity proof of Lazard's
  method for CAD construction}. J. Symb. Comput. 92, 52--69.

\bibitem[{Monagan(2004)}]{M04}
Monagan, M., 2004. {Maximal Quotient Rational Reconstruction: An Almost Optimal
  Algorithm for Rational Reconstruction}. In: Proc. {ISSAC, ACM 2004}. pp.
  243--249.

\bibitem[{Mora and Robbiano(1993)}]{MR88}
Mora, T., Robbiano, L., 1993. {The Gr\"obner Fan of an Ideal}. J. Symb. Comput.
  15, 199--209.

\bibitem[{Noro(2002)}]{Nor2002}
Noro, M., 2002. {An efficient modular algorithm for computing the global {\rm
  b}-function}. In: Proceeding of the First Congress of Mathematical Software,
  A.M.~Cohen, X.S.~Gao, N.~Takayama, eds. pp. 147--157.

\bibitem[{Noro and Yokoyama(1999)}]{NY[17]}
Noro, M., Yokoyama, K., 1999. {A modular method to compute the rational
  univariate representation of zero-dimensional ideals}. J. Symb. Comput. 28,
  243--263.

\bibitem[{Noro and Yokoyama(2004)}]{NY[18]}
Noro, M., Yokoyama, K., 2004. {Implementation of prime decomposition of
  polynomial ideals over small finite fields}. J. Symb. Comput. 38, 1227--1246.

\bibitem[{Noro and Yokoyama(2018)}]{NY2018}
Noro, M., Yokoyama, K., 2018. {Usage of Modular Techniques for Efficient
  Computation of Ideal Operations}. Math.Comput.Sci. 12, 1--32.

\bibitem[{Pauer(2007)}]{Pa}
Pauer, F., 2007. {Gr\"obner bases with coefficients in rings}. J. Symb. Comput.
  42, 1003--1011.

\bibitem[{Rouillier(1999)}]{Rou}
Rouillier, F., 1999. {Solving Zero-Dimensional Systems Through the Rational
  Univariate Representation}. AAECC 9, 433--461.

\bibitem[{Winkler(1988)}]{Wi}
Winkler, F., 1988. {A $p$-adic Approach to the Computation of Gr\"obner Bases}.
  J. Symb. Comput. 6, 287--304.

\end{thebibliography}

\end{document}